\documentclass[onefignum,onetabnum]{siamart190516}
\usepackage{amsmath}
\usepackage{amssymb}
\usepackage{enumerate}
\usepackage{mathrsfs}
\usepackage{cases}
\usepackage{color}
\usepackage{tikz}
\usepackage{makecell}
\usepackage{multirow}
\usepackage{siunitx}
\usepackage{algorithmic,algorithm}
\usepackage{graphicx}
\usepackage{mathtools}
\usepackage{pgfplots}
\usepackage{scalefnt}
\usepgfplotslibrary{fillbetween}
\definecolor{darkblue}{rgb}{0.0, 0.0, 0.55}
\definecolor{bordeaux}{rgb}{0.34, 0.01, 0.1}
\definecolor{color1}{RGB}{145,30,180}
\definecolor{color2}{RGB}{245,130,48}
\definecolor{color3}{RGB}{230,25,75}
\definecolor{lightgreen}{RGB}{144, 238, 144}
\newsiamthm{example}{Example}

\def\Z{{\mathbb{Z}}}

\def\R{{\mathbb{R}}}
\def\C{{\mathbb{C}}}
\def\N{{\mathbb{N}}}

\def\x{{\mathbf{x}}}

\def\S{{\mathbb{S}}}

\def\Tr{\hbox{\rm{Tr}}}

\def\grad{\hbox{\rm{grad}}}
\def\hess{\hbox{\rm{Hess}}}

\def\rank{\hbox{\rm{rank}}}

\def\cA{{\mathcal{A}}}
\def\cI{{\mathcal{I}}}
\def\cB{{\mathcal{B}}}
\def\cM{{\mathcal{M}}}
\def\cN{{\mathcal{N}}}
\def\Diag{\hbox{\rm{Diag}}}

\DeclareMathOperator*{\argmin}{arg\,min}

\usepackage{lipsum}
\usepackage{amsfonts}
\usepackage{graphicx}
\usepackage{epstopdf}
\usepackage{algorithmic}
\ifpdf
  \DeclareGraphicsExtensions{.eps,.pdf,.png,.jpg}
\else
  \DeclareGraphicsExtensions{.eps}
\fi

\newsiamremark{remark}{Remark}
\newsiamremark{hypothesis}{Hypothesis}
\crefname{hypothesis}{Hypothesis}{Hypotheses}
\newsiamthm{claim}{Claim}
\newsiamthm{assumption}{Assumption}

\headers{A Dual Riemannian ALM for Low-Rank SDPs with Unit Diagonal}{J. Wang, L. Hu, and B. Xia}

\title{A Dual Riemannian Augmented Lagrangian Method for Low-Rank SDPs with Unit Diagonal
\thanks{Submitted to the editors DATE.
}}

\author{Jie Wang\thanks{State Key Laboratory of Mathematical Sciences, Academy of Mathematics and Systems Science, Chinese Academy of Sciences, Beijing, China
  (\email{wangjie212@amss.ac.cn})} \and Liangbing Hu\thanks{Nanjing Research Institute of Electronics Technology, Nanjing, China
  (\email{huliangbing2000@163.com})}
\and Bican Xia\thanks{(School of Mathematical Sciences, Peking University, Beijing, China (\email{xbc@math.pku.edu.cn})}
  }

\usepackage{amsopn}
\DeclareMathOperator{\diag}{diag}

\ifpdf
\hypersetup{
  pdftitle={A Dual Riemannian ALM for Low-Rank SDPs with Unit Diagonal},
  pdfauthor={J. Wang, L. Hu, and B. Xia}
}
\fi

\begin{document}

\maketitle

\begin{abstract}
We propose \texttt{ManiDSDP}, a dual Riemannian augmented Lagrangian method for solving low-rank semidefinite programs in the dual form whose positive semidefinite variable has unit diagonal. The method first eliminates the unconstrained variable by a variable projection, expresses affine consistency through an orthogonal-projector residual, and then applies the Burer--Monteiro factorization. The resulting subproblems are optimization problems on the oblique manifold and are solved by a Riemannian trust-region method with dynamic rank adjustment. A negative-curvature direction extracted from the certificate of a completed subproblem is carried into the next subproblem, leading to a one-step-delayed curvature-correction design. 
Under controlled inner errors, we prove that every cluster point returned by the method is a KKT point. Under additional full-sequence convergence, strict complementarity, and finite-tail exact complementarity, the iterations identify the rank of the limiting solution in finitely many steps. Under a strictly complementary rank-one solution, we further establish a local one-step locking theorem, providing a mechanism for the observed residue diving phenomenon. Numerical experiments on dense and sparse BQPs and on UCQPs demonstrate high accuracy and favorable efficiency and scalability relative to \texttt{MOSEK}, \texttt{SDPNAL+}, and \texttt{ManiSDP}.
\end{abstract}



\section{Introduction}
In this paper, we consider the semidefinite program (SDP)
\begin{equation}\label{dsdp}
\begin{cases}
\inf\limits_{y\in\R^m,S\in\S_n}&b^{\intercal}y\\
\,\,\,\,\quad\rm{s.t.}&S=\cA^*(y)-C\succeq0,\\
&\diag(S)=\mathbf 1,
\end{cases}\tag{DSDP}
\end{equation}
where $\cA:\S_n\to\R^m$ is a linear operator, $\cA^*$ is its adjoint, and $\S_n$ denotes the space of $n\times n$ real symmetric matrices. We focus on the following structured regime.
\begin{assumption}\label{assump1}
Problem~\eqref{dsdp} admits a low-rank optimal solution $S^{\star}$.
\end{assumption}
\begin{assumption}\label{assump2}
The linear operator $\cA\cA^*$ is invertible.
\end{assumption}
Assumption~\ref{assump1} is a structural assumption on the SDP. In particular, moment relaxations in the dual form (which fit in \eqref{dsdp}) of polynomial optimization problems are observed or expected to satisfy Assumption~\ref{assump1}.  Assumption~\ref{assump2} is equivalent to injectivity of $\cA^*$, or, equivalently, to linear independence of the coefficient matrices defining $\cA$. Redundant affine constraints can therefore be removed in preprocessing.

Problems of the form~\eqref{dsdp} arise as moment-sum-of-squares (SOS) relaxations of real, complex, and noncommutative polynomial optimization problems, including binary quadratic programs (BQPs) \cite{wang2025solving}, the classical Ising problem \cite{edwards1975theory}, the generalized orthogonal Procrustes problem \cite{ling2025generalized}, group synchronization \cite{ling2023solving,ling2025local,mcrae2025benign}, inference in graphical models \cite{erdogdu2017inference}, synchronization and community detection \cite{bandeira2016low}, unit-modulus complex quadratic programs (UCQPs) \cite{soltanalian2014designing,waldspurger2015phase}, multiple-input multiple-output detection \cite{lu2019tightness}, and ground-state problems for quantum spin systems \cite{wang2023certifying}.

Despite convexity, solving large SDPs to high accuracy remains difficult. Interior-point methods are accurate and robust on small and medium instances, but their memory and arithmetic costs become prohibitive at scale \cite{andersen2003implementing,toh1999sdpt3}. First-order splitting methods, including alternating direction methods of multipliers (ADMM), have lower per-iteration costs but may converge slowly and often return only low- or medium-accuracy solutions \cite{garstka2021cosmo,wen2010alternating,zheng2020chordal}. Augmented Lagrangian methods (ALM) with semismooth Newton or related inner solvers offer a stronger accuracy--scalability tradeoff \cite{yang2015sdpnal,zhao2010newton}. On the other hand,
SDPs from practice frequently possess certain structures, e.g., low-rank optimal solutions, data sparsity, and unit diagonal constraints. Those structures could be exploited
to design more efficient and scalable algorithms \cite{hou2025low,hou2025rinnal+,lemon2016low,li2026curvature,yang2022inexact}.

The low-rank property is commonly exploited through the Burer--Monteiro factorization \cite{burer2003nonlinear,burer2005local}; for SDPs whose constraints define smooth manifolds, this factorization leads naturally to Riemannian optimization methods \cite{journee2010low,wang2025solving,wang2023decomposition}. In particular, \cite{wang2025solving} introduced \texttt{ManiSDP}, a primal Riemannian ALM with dynamic rank adjustment and negative-curvature directions for escaping nonoptimal stationary points. The present work adopts those two ingredients, but differs in two essential respects: it works on the dual side and performs a variable projection before factorization. These choices lead to different subproblems, multiplier invariants, and convergence conditions from those of primal low-rank ALM methods and direct full-matrix ADMM splittings.

There are two principal reasons for working on the dual side. First, the cost of solving an SDP depends not only on the positive-semidefinite (PSD) block size but also on the number of independent affine constraints, equivalently the rank of $\cA$. In higher-order moment--SOS relaxations, the dual SOS formulation can have substantially fewer affine constraints than the primal moment formulation. Section~\ref{sec3-1} quantifies this comparison for second-order relaxations of BQPs and UCQPs.
Second, Assumption~\ref{assump2} permits a variable projection: $y(S)\coloneqq(\cA\cA^*)^{-1}\cA(S+C)$.
It yields an equivalent problem in $S$ with the projected residual
$(\cI-\cA^*(\cA\cA^*)^{-1}\cA)(S+C)=0$.
We show that applying the two-block ADMM directly to~\eqref{dsdp} produces, after the $y$-update, the same genuine feasibility penalty plus an additional proximal term in $\operatorname{range}(\cA^*)$. This term damps changes in directions that affect the objective but can already be represented exactly by $y$. The variable-projected formulation instead updates $y(S)$ exactly for every candidate $S$ and penalizes only genuine affine inconsistency.

Our contributions are summarized as follows.

$\bullet$ We develop \texttt{ManiDSDP}, a dual Riemannian augmented Lagrangian method for solving \eqref{dsdp}. After the variable projection and the factorization $S=YY^{\intercal}$ with $Y\in\R^{n\times p}$, the unit-diagonal condition becomes the oblique-manifold constraint. A negative eigenvalue of the certificate of the completed factorized subproblem yields an explicit second-order descent direction after increasing the factorization size. The \texttt{ManiDSDP} method carries this direction into the next outer iteration. We derive the exact certificate perturbation caused by this one-step delay and give a computable condition under which the carried direction remains a negative-curvature direction for the new subproblem.

$\bullet$ We establish global convergence of \texttt{ManiDSDP} under a one-step-delayed forcing condition. A Lyapunov argument then proves summability of the negative-eigenvalue defects, boundedness of all primal-dual sequences, convergence of the entire multiplier sequence, and convergence of the whole sequence to the KKT set. A summable spectral-tolerance condition is obtained as a direct corollary.

$\bullet$ We prove finite rank identification under full-sequence convergence, strict complementarity, and an exact finite-tail complementarity condition. The latter holds, for example, when the Riemannian subproblems are solved exactly to first-order criticality on a tail of the sequence or when dual feasibility is eventually exact.

$\bullet$ Under strict complementarity and a rank-one optimal solution, we establish a local one-step locking theorem. It shows that once the multiplier enters an explicit neighborhood, every exact subproblem has the rank-one solution as its unique global minimizer, independently of the penalty parameter, and the subsequent iterations remain locked there. This provides a theoretical mechanism for the residue diving phenomenon observed in the numerical experiments.

$\bullet$ Various numerical experiments are performed to evaluate the performance of the proposed algorithm. We demonstrate that it outperforms, by a significant margin on the tested medium- and large-scale instances, two advanced SDP solvers (\texttt{MOSEK}, \texttt{SDPNAL+}) in terms of accuracy, efficiency, and scalability on second-order SDP relaxations of dense and sparse BQPs and UCQPs. In addition, it runs several times faster than the recent low-rank SDP solver {\tt ManiSDP} \cite{wang2025solving}. The experiments also indicate an intriguing feature of our algorithm applied to random BQPs — the residue diving phenomenon, that is, the maximal KKT residue sharply decreases to far less than $10^{-8}$ at some ALM iteration. Moreover, it is worth mentioning that our algorithm typically returns a solution of extremely high accuracy in a few tens of ALM iterations. In comparison, traditional ADMM-based algorithms often require more than ten thousand iterations to arrive at a solution of high accuracy for the same type of SDPs \cite{kang2025local}.

Several other recent low-rank SDP methods are particularly relevant to the present work. The Burer--Monteiro ADMM of Chen and Goulart \cite{chen2022burer} and \texttt{LoRADS} \cite{han2025low} introduce bilinear matrix splittings so that the ADMM subproblems are quadratic: the former is tailored to diagonal constraints, whereas the latter treats general SDPs. \texttt{SDPLR+} \cite{huang2024suboptimality} and \texttt{HALLaR} \cite{monteiro2024low} instead operate on primal trace-bounded SDPs: \texttt{SDPLR+} uses a computable suboptimality bound to control termination and rank growth, while \texttt{HALLaR} combines an inexact ALM with accelerated proximal-point and Frank--Wolfe steps. The feasible Riemannian method of Tang and Toh \cite{tang2024feasible} maintains affine feasibility and adapts both rank and support without a penalty parameter, whereas \texttt{RiNNAL-POP} \cite{hou2026low} alternates between low-rank and convex-lifting phases and exploits facial, polyhedral, and moment-consistency structures in moment-SOS relaxations. By comparison, \texttt{ManiDSDP} is specialized to the dual unit-diagonal model~\eqref{dsdp}. This specialization is less general---but it can reduce the affine-constraint dimension and permits the specific multiplier and KKT-set convergence analysis, finite rank identification under strict complementarity, and local one-step locking.

The rest of the paper is organized as follows. Section~\ref{preliminaries} introduces notation and preliminary results. Section~\ref{ALM} develops the variable-projected dual augmented Lagrangian framework, compares the affine dimensions of primal moment and dual SOS formulations, and explains the advantage of the projected formulation over the direct formulation. Section~\ref{manifold} derives the Riemannian geometry, gradient, Hessian, and negative-curvature directions for the factorized subproblems. Section~\ref{algorithms} presents \texttt{ManiDSDP}, analyzes the one-step-delayed use of the curvature direction, and establishes global convergence and finite rank identification. Section~\ref{finite-locking} proves the local one-step locking theorem and the inexact near-locking bounds. Section~\ref{experiments} reports numerical results, and Section~\ref{conclusions} concludes the paper.

\section{Notation and preliminaries}\label{preliminaries}
We write $\R$, $\R_+$, and $\N$ for the real numbers, the nonnegative real numbers, and the nonnegative integers, respectively. The spaces of $n\times n$ real symmetric and PSD matrices are denoted by $\S_n$ and $\S_n^+$, respectively. For a matrix $A$, $A^{\intercal}$ denotes its transpose; for square $A$, $\Tr(A)$ denotes its trace. For two matrices $A,B\in\R^{n\times m}$, the inner product is defined as $\langle A, B\rangle=\Tr(A^{\intercal}B)$. 
For $A\in\R^{n\times n}$, $\diag(A)$ denotes the diagonal of $A$, $\Diag(A)$ denotes the diagonal matrix with the same diagonal as $A$, and $\lambda_{\min}(A),\lambda_{\max}(A)$ denote the smallest, largest eigenvalues of $A$, respectively. For a vector $v$, let $|v|$ be the dimension of $v$, $\|v\|$ the $2$-norm of $v$, and $\Diag(v)$ the diagonal matrix with $v$ being its diagonal. For a matrix $A$, $\|A\|$ (resp. $\|A\|_2$) denotes the Frobenius (resp. spectral) norm of $A$.
For a function $f(x)$, we write $\nabla f(x)$ (resp. $\grad\,f(x)$) for the Euclidean (resp. Riemannian) gradient, and write $\nabla^2 f(x)[u]$ (resp. $\hess\,f(x)[u]$) for the Euclidean (resp. Riemannian) Hessian acting on $u$. For a set $E$, $|E|$ denotes its cardinality.

Let $\mathcal{A}(X)\coloneqq(\langle A_i, X\rangle)_{i=1}^m$ for $A_i\in\S_n$ and $\cA^*:\R^m\rightarrow\S_n$ be the adjoint operator of $\cA$ defined as $\cA^*(y)\coloneqq\sum_{i=1}^my_iA_i$ for $y\in\R^m$. The linear operator $\cA\cA^*:\R^m\rightarrow\R^m$ is defined as $\cA\cA^*(y)=\cA(\cA^*(y))$. Let $\mathcal{B}(S)\coloneqq(\langle B_i, S\rangle)_{i=1}^n$, where $B_i\in\S_n$ is the matrix with the $i$-th diagonal element being one and all other elements being zeros, and $\cB^*:\R^n\rightarrow\S_n$ be the adjoint operator of $\cB$ defined as $\cB^*(z)\coloneqq\sum_{i=1}^nz_iB_i=\Diag(z)$ for $z\in\R^n$. We use $\cI$ to stand for the identity operator.

The dual of \eqref{dsdp} reads as
\begin{equation}\label{psdp}
\begin{cases}
\sup\limits_{X\in\S_n^+,z\in\R^n}&\langle C, X+\cB^*(z)\rangle+\sum_{i=1}^nz_i\\
\,\,\,\,\,\,\,\,\,\rm{s.t.}&\mathcal{A}(X+\cB^*(z))=b.
\end{cases}\tag{PSDP}
\end{equation}
Throughout the paper, we assume that strong duality holds for \eqref{psdp}--\eqref{dsdp}. A quadruple $(S,y,X,z)$ is a KKT point of \eqref{dsdp} if
\begin{subequations}\label{eq:kkt-system}
\begin{align}
    S&=\cA^*(y)-C\succeq0, & \diag(S)&=\mathbf1,\\
    X&\succeq0, & \cA(X+\cB^*(z))&=b,\\
    XS&=0.&&
\end{align}
\end{subequations}
We denote the set of such quadruples by $\mathcal K$.


\section{The dual projected ALM approach}\label{ALM}
In this section, we introduce the dual projected ALM framework for solving \eqref{dsdp}.
Let us define
\begin{equation}\label{maniX}
\cM\coloneqq\{S\in\S_n^+\mid\diag(S)=\mathbf{1}\}=\{S\in\S_n^+\mid\cB(S)=\mathbf{1}\}.
\end{equation}
Throughout the paper, we define the constant matrix $D$ as
\begin{equation}
D\coloneqq\cA^*(\cA\cA^*)^{-1}(b).
\end{equation}
From $S=\cA^*(y)-C$ together with Assumption~\ref{assump2}, we have\footnote{Indeed, $y(S)=(\cA\cA^*)^{-1}\cA(S+C)$ is the unique solution to the least-squares problem
$\min_{y}\frac{1}{2}\|\mathcal{A}^*(y)-S-C\|^2$.}
\[y=y(S)\coloneqq(\cA\cA^*)^{-1}\cA(S+C),\]
and consequently, 
\begin{equation*}
b^{\intercal}y=b^{\intercal}(\cA\cA^*)^{-1}\cA(S+C)=\left\langle\cA^*(\cA\cA^*)^{-1}(b), S+C\right\rangle=\langle D, S+C\rangle.    
\end{equation*}
Therefore, we may reformulate \eqref{dsdp} as
\begin{equation}\label{dsdp1}
\begin{cases}
\inf\limits_{S\in\cM}&\langle D, S+C\rangle\\
\,\,\,\rm{s.t.}&(\cI-\cA^*(\cA\cA^*)^{-1}\cA)(S+C)=0.
\end{cases}\tag{DSDP'}
\end{equation}
From now on, we choose to work on the reformulation \eqref{dsdp1} instead of \eqref{dsdp}. We employ the ALM framework to solve \eqref{dsdp1}.
The augmented Lagrangian function associated with \eqref{dsdp1} is defined by
\begin{equation}
   L_{\sigma}(S,y,\widetilde{X})=\langle D, S+C\rangle-\langle \widetilde{X}, \cA^*(y)-S-C\rangle+\frac{\sigma}{2}\| \cA^*(y)-S-C\|^2.
\end{equation}
At the $k$-th ALM iteration, one solves the subproblem:
\begin{equation}
S^{k+1}\coloneqq\argmin\limits_{S\in\cM}L_{\sigma_k}\left(S, y(S), \widetilde{X}^k\right).\label{subp1}
\end{equation}
Then, $y^{k}$ and $\widetilde{X}^k$ are updated as
\begin{subequations}
\begin{align}
y^{k+1}&\coloneqq(\cA\cA^*)^{-1}\cA(S^{k+1}+C),\\
\widetilde{X}^{k+1}&\coloneqq \widetilde{X}^k-\sigma_k\left(\cA^*\left(y^{k+1}\right)-S^{k+1}-C\right).
\end{align}
\end{subequations}

Let $\Phi_k(S)\coloneqq L_{\sigma_{k}}(S,y(S), \widetilde{X}^{k})$.
The following lemma characterizes the optimality conditions of \eqref{subp1}.
\begin{lemma}\label{lm:suboptimality}
An element $S\in\cM$ is a minimizer of \eqref{subp1} if and only if
\begin{subequations}
\begin{align}
X\coloneqq \nabla\Phi_k(S)-\Diag(\nabla\Phi_k(S)S)&\succeq0,\label{optimality3}\\ 
XS&=0, \label{optimality4}
\end{align}
\end{subequations}
where $\nabla\Phi_k(S)=\widetilde{X}^k+D-\sigma_k(\cA^*(y(S))-S-C)$.
\end{lemma}
\begin{proof}
Since \eqref{subp1} is convex, $S\in\cM$ is a minimizer if and only if the KKT conditions hold.
The Lagrangian function of \eqref{subp1} is
\begin{equation}
    \Phi_k(S)-z^{\intercal}(\cB(S)-\mathbf{1})-\langle X,S\rangle,
\end{equation}
where $z\in\R^n$ and $X\in\S_n^+$ are the dual variables.
The KKT conditions of \eqref{subp1} thus read as
\begin{subequations}
\begin{align}
S&\in\cM,\\
\nabla\Phi_k(S)-\cB^*(z)-X&=0,\\
X&\succeq0,\\
XS&=0.
\end{align}
\end{subequations}
Substituting $X=\nabla\Phi_k(S)-\cB^*(z)=\nabla\Phi_k(S)-\Diag(z)$ into $XS=0$ yields
\begin{equation}\label{sec3:eq1}
\nabla\Phi_k(S)S=\Diag(z)S.
\end{equation}
Taking diagonal entries in \eqref{sec3:eq1} and using $\diag(S)=\mathbf1$ gives
$z=\diag(\nabla\Phi_k(S)S)$.
\end{proof}

According to Lemma \ref{lm:suboptimality}, 
$X^k$ is updated as
\begin{equation}
X^{k+1}\coloneqq\widetilde{X}^{k+1}+D-\Diag\left(\left(\widetilde{X}^{k+1}+D\right)S^{k+1}\right).
\end{equation}
We present the ALM algorithm for solving \eqref{dsdp1} in Algorithm \ref{alg0}.

\begin{algorithm}
\renewcommand{\algorithmicrequire}{\textbf{Input:}}
\renewcommand{\algorithmicensure}{\textbf{Output:}}
\caption{The ALM algorithm for solving \eqref{dsdp1}}\label{alg0}
	\begin{algorithmic}[1]
   	\REQUIRE $\cA,b,C$; $\sigma_0>0$; $0<\sigma_{\min}<\sigma_{\max}$
		\STATE $k \gets 0$, $\widetilde{X}^0\gets 0_{n\times n}$;
		\WHILE{The stopping criterion is not fulfilled}
		\STATE Obtain $S^{k+1}$ by solving the ALM subproblem \eqref{subp1};
		\STATE $y^{k+1} \gets (\cA\cA^*)^{-1}\cA(S^{k+1}+C)$;
            \STATE $\widetilde{X}^{k+1}\gets \widetilde{X}^k-\sigma_k\left(\cA^*\left(y^{k+1}\right)-S^{k+1}-C\right)$;
            \STATE $z^{k+1}\gets\diag\left(\left(\widetilde{X}^{k+1}+D\right)S^{k+1}\right)$;
            \STATE $X^{k+1}\gets\widetilde{X}^{k+1}+D-\Diag\left(z^{k+1}\right)$;
            \STATE Determine $\sigma_{k+1}\in[\sigma_{\min},\sigma_{\max}]$ according to certain policy;
		\STATE $k \gets k+1$;
		\ENDWHILE
		\RETURN $(S^{k},y^{k},X^{k},z^{k})$;
	\end{algorithmic}
\end{algorithm}


We prove the following two lemmas regarding Algorithm~\ref{alg0}.

\begin{lemma}\label{lm1}
For every $k\ge0$, $\cA(\widetilde{X}^{k})=0$. Furthermore, $\cA(X^{k}+\Diag(z^k))=b$ for every $k\geq1$.
\end{lemma}
\begin{proof}
By the definition of $y^{k+1}$ and the invertibility of $\cA\cA^*$,
\begin{align*}
\cA\left(\cA^*(y^{k+1})-S^{k+1}-C\right)
&=(\cA\cA^*)(y^{k+1})-\cA(S^{k+1}+C)\\
&=(\cA\cA^*)(\cA\cA^*)^{-1}\cA(S^{k+1}+C)
  -\cA(S^{k+1}+C)\\
&=0.
\end{align*}
Applying $\cA$ to the multiplier update therefore gives
\begin{equation*}
\cA(\widetilde{X}^{k+1})=\cA(\widetilde{X}^{k})
 -\sigma_k\cA\left(\cA^*(y^{k+1})-S^{k+1}-C\right)=\cA(\widetilde{X}^{k}).
\end{equation*}
Hence $\cA(\widetilde{X}^{k})=\cA(\widetilde{X}^{0})=0$ for every $k\ge0$.
Furthermore,
\begin{equation*}
    \cA(X^{k}+\Diag(z^k))=\cA(\widetilde{X}^{k}+D)=\cA(D)=b.
\end{equation*}
\end{proof}

\begin{lemma}\label{sec3:lm1}
For every $k\ge1$, $\langle X^{k}, S^{k}\rangle=0$.
\end{lemma}
\begin{proof}
We have
\begin{align*}
\left\langle X^{k}, S^{k} \right\rangle
&=
\left\langle
\widetilde{X}^{k}+D-\operatorname{Diag}\left(z^{k}\right),
S^{k}\right\rangle\\
&=
\operatorname{Tr}\left(
\left(\widetilde{X}^{k}+D\right)S^{k}
\right)
-
\left(z^{k}\right)^{\intercal}
\operatorname{diag}\left(S^{k}\right)
\\
&=
\mathbf{1}^{\intercal}
\operatorname{diag}\left(
\left(\widetilde{X}^{k}+D\right)S^{k}
\right)
-
\mathbf{1}^{\intercal}z^{k}
\\
&=0,
\end{align*}
where we have used $\operatorname{diag}\left(S^{k}\right)=\mathbf{1}$.
\end{proof}

\subsection{Why solve the dual formulation}\label{sec3-1}
The computational cost of an SDP depends not only on the size of PSD blocks but also on the number of independent affine constraints, equivalently the rank of $\cA$. We next compare this quantity for the second-order dual SOS and primal moment relaxations of two representative classes of polynomial optimization problems to illustrate the advantage of the dual formulation.

\paragraph{Binary quadratic programs}
Consider the BQP 
\begin{equation}\label{sec3:bqp}
\inf_{\x\in\R^q}\ \  \x^\intercal Q\x + \mathbf{c}^\intercal\x
\qquad\mathrm{s.t.}\qquad x_i^2=1,\quad i=1,\ldots,q,\tag{BQP}
\end{equation}
where $Q\in\S_q$ and $\mathbf{c}\in\R^{q}$. Let
\begin{equation*}
v(\x)\coloneqq[1,x_1,\ldots,x_q,x_1x_2,x_1x_3,\ldots,x_{q-1}x_q]^{\intercal}
\end{equation*}
be the vector of monomials in $\x$ up to degree two (excluding squares $x_i^2,\,i=1,\ldots,q$). The second-order dual SOS relaxation of \eqref{sec3:bqp} is the following SDP:
\begin{equation}\label{sec6:sos}
\begin{cases}
\sup\limits_{X,\lambda} &\lambda\\
\,\rm{s.t.}&\x^\intercal Q\x + \mathbf{c}^\intercal\x-\lambda\equiv v(\x)^{\intercal}Xv(\x)\quad\mathrm{mod}\,\,I,\\
&X\succeq0,\\
\end{cases}
\end{equation}
where $I$ is the ideal generated by $\{x_i^2-1\}_{i=1}^q$ in the polynomial ring $\R[\x]$. Each linear constraint of \eqref{sec6:sos} is indexed by a monomial in \[\{1\}\cup\{x_i\}_{i=1}^q\cup\{x_ix_j\}_{1\le i<j\le q}\cup\{x_ix_jx_k\}_{1\le i<j<k\le q}\cup\{x_ix_jx_kx_l\}_{1\le i<j<k<l\le q}.\] 
Hence, in the dual formulation of \eqref{sec6:sos}, the rank of $\cA$ is
\begin{equation}
    m^{\rm B}_{\rm sos}
    =
    \frac{q^4-2q^3+11q^2+14q+24}{24}.
    \label{eq:bqp-number-moments}
\end{equation}
Therefore, in the primal moment SDP, the rank of the corresponding $\cA$ is
\begin{equation}
    m^{\rm B}_{\rm mom}
    =\frac{q^4+4q^3+5q^2+2q+12}{12}.
    \label{eq:bqp-moment-number-constraints}
\end{equation}
In particular,
\begin{equation}
    m^{\rm B}_{\rm sos}\sim \frac{q^4}{24},
    \qquad
    m^{\rm B}_{\rm mom}\sim \frac{q^4}{12},
    \qquad
    \frac{m^{\rm B}_{\rm mom}}{m^{\rm B}_{\rm sos}}\longrightarrow 2.
    \label{eq:bqp-size-ratio}
\end{equation}
Thus, the dual SOS formulation asymptotically halves the rank of $\cA$.

\paragraph{Unit-modulus complex quadratic programs}
Consider the real-coefficient UCQP
\begin{equation}
    \inf_{\x\in\C^q}\ \ \x^*Q\x + \operatorname{Re}(\mathbf{c}^\intercal\x)
    \qquad
    \mathrm{s.t.}\qquad
    |x_i|^2=1,\quad i=1,\ldots,q,\tag{UCQP}
    \label{eq:ucqp-size-comparison}
\end{equation}
where $Q\in\S_q$ and $\mathbf{c}\in\R^{q}$. Because the coefficients are real, the moment-SOS relaxation for \eqref{eq:ucqp-size-comparison} can be taken to be real \cite{wang2025real}. 
We shall use the reduced monomial basis $\tilde{v}(\x)\coloneqq[\x,\bar \x]_2\setminus\{x_i\bar x_i:i\in[q]\}$\footnote{$[\x,\bar \x]_2$ denotes the set of monomials in $\x$ and $\bar \x$ up to degree two.} with $|\tilde{v}(\x)|=\sum_{j=0}^2 2^j\binom qj\binom2j=2q^2+2q+1$ to build the moment-SOS relaxation.
Products of two reduced basis monomials generate precisely the Laurent exponent set
$\Delta_4=\left\{\delta\in\Z^q:\|\delta\|_1\leq4\right\}$
with
\begin{equation}
|\Delta_4|=\sum_{j=0}^{4}2^j\binom qj\binom4j=\frac{2q^4+4q^3+10q^2+8q+3}{3}.
\end{equation}
Because the SDP is real, the moments and coefficients indexed by $\delta$ and $-\delta$ coincide. The number of independent real Laurent moments, or equivalently independent real coefficient classes, is thus
\begin{equation}
    m^{\rm U}_{\rm sos}
    =
    \frac{|\Delta_4|+1}{2}=\frac{q^4+2q^3+5q^2+4q+3}{3}.
    \label{eq:ucqp-reduced-number-real-coefficients}
\end{equation}
Hence, the primal moment SDP has
\begin{equation}
    m^{\rm U}_{\rm mom}
    =\frac{5q^4+10q^3+10q^2+5q+3}{3}
    \label{eq:ucqp-reduced-moment-number-constraints}
\end{equation}
independent real affine constraints.
Consequently,
\begin{equation}
    m^{\rm U}_{\rm sos}\sim\frac{q^4}{3},
    \qquad
    m^{\rm U}_{\rm mom}\sim\frac{5q^4}{3},
    \qquad
    \frac{m^{\rm U}_{\rm mom}}{m^{\rm U}_{\rm sos}}
    \longrightarrow5.
    \label{eq:ucqp-reduced-size-ratio}
\end{equation}
Thus, the dual SOS formulation has asymptotically one fifth as many independent affine constraints as the primal moment formulation. 


\subsection{Why solve \eqref{dsdp1} instead of \eqref{dsdp}}
Although \eqref{dsdp1} and \eqref{dsdp} are equivalent optimization problems, applying the ALM algorithm or the ADMM algorithm to the two formulations leads to different subproblems and consequently different numerical behaviors. In this subsection, we show that applying ADMM directly to \eqref{dsdp} introduces an additional proximal regularization term that is absent from the subproblem associated with \eqref{dsdp1}. This term penalizes changes in directions that do not correspond to constraint violation and may therefore substantially slow down the algorithm.

Define linear operators
\begin{equation}
    \mathcal{P}
    \coloneqq
    \mathcal{A}^{*}
    (\mathcal{A}\mathcal{A}^{*})^{-1}
    \mathcal{A},
    \qquad
    \mathcal{Q}
    \coloneqq
    \mathcal{I}-\mathcal{P}.
    \label{eq:projection-operators}
\end{equation}
Under Assumption~\ref{assump2}, $\mathcal{P}$ is the orthogonal projector onto $\operatorname{range}(\mathcal{A}^{*})$, whereas $\mathcal{Q}$ is the orthogonal projector onto $\ker(\mathcal{A})$. Let $Z\coloneqq S+C$.
Then
\begin{equation}
    \mathcal{A}^{*}(y(S))
    =
    \mathcal{P}(Z),
    \qquad
    \mathcal{A}^{*}(y(S))-Z
    =
    -\mathcal{Q}(Z).
    \label{eq:projected-y}
\end{equation}
In particular,
\begin{equation}
    \min_{y\in\mathbb{R}^{m}}
    \left\|
        \mathcal{A}^{*}(y)-Z
    \right\|^{2}
    =
    \left\|\mathcal{Q}(Z)\right\|^{2}.
    \label{eq:projection-distance}
\end{equation}
Thus, $\|\mathcal{Q}(Z)\|$ is precisely the distance from $Z$ to
$\operatorname{range}(\mathcal{A}^{*})$ and measures the genuine violation of the
constraint $Z\in\operatorname{range}(\mathcal{A}^{*})$.

By Lemma~\ref{lm1}, $\mathcal{A}(\widetilde X^{k})=0$, and hence
$\widetilde X^{k}\in\ker(\mathcal{A})$. Moreover,
$D\in\operatorname{range}(\mathcal{A}^{*})$. Using
\eqref{eq:projected-y}, the objective of the subproblem~\eqref{subp1} can therefore be written as
\begin{align}
    \Phi_k(S)
    &=
    \left\langle D,Z\right\rangle
    -
    \left\langle
        \widetilde X^{k},
        \mathcal{P}(Z)-Z
    \right\rangle
    +
    \frac{\sigma_k}{2}
    \left\|\mathcal{P}(Z)-Z\right\|^{2}
    \notag\\
    &=
    \left\langle
        D+\widetilde X^{k},
        Z
    \right\rangle
    +
    \frac{\sigma_k}{2}
    \left\|\mathcal{Q}(Z)\right\|^{2}.
    \label{eq:reduced-S-subproblem}
\end{align}
Consequently, at the $k$-th ALM iteration for \eqref{dsdp1}, one solves
\begin{equation}
    \min_{S\in\mathcal{M}}
    \left\{
        \left\langle
            D+\widetilde X^{k},
            S+C
        \right\rangle
        +
        \frac{\sigma_k}{2}
        \left\|
            \mathcal{Q}(S+C)
        \right\|^{2}
    \right\}.
    \label{eq:reduced-ALM-subproblem}
\end{equation}

Next, we compare \eqref{eq:reduced-ALM-subproblem} with the subproblem obtained by applying the standard two-block ADMM directly to \eqref{dsdp}. The augmented
Lagrangian for the direct formulation is
\begin{equation}
    \widehat L_{\sigma}(S,y,\overline X)
    \coloneqq
    b^{\intercal}y
    -
    \left\langle
        \overline X,
        \mathcal{A}^{*}(y)-S-C
    \right\rangle
    +
    \frac{\sigma}{2}
    \left\|
        \mathcal{A}^{*}(y)-S-C
    \right\|^{2}.
    \label{eq:direct-augmented-lagrangian}
\end{equation}
Let $\overline X^{k}$ denote the multiplier generated by this direct ADMM scheme. The first-order optimality condition of its $y$-subproblem, together with
the multiplier update, gives
\[
    \mathcal{A}(\overline X^{k+1})=b.
\]
Hence, after the first multiplier update, one may write
\begin{equation}
    \overline X^{k}=D+W^{k},
    \qquad
    \mathcal{A}(W^{k})=0.
    \label{eq:shifted-direct-multiplier}
\end{equation}
At every subsequent iteration, the optimality condition of the $y$-subproblem
also yields
\begin{equation}
    \mathcal{A}^{*}(y^{k})
    =
    \mathcal{P}(Z^{k}),
    \qquad
    Z^{k}\coloneqq S^{k}+C.
    \label{eq:direct-y-projection}
\end{equation}
Indeed, once $\mathcal{A}(\overline X^{k-1})=b$, the $y$-optimality condition
reduces to $(\mathcal{A}\mathcal{A}^{*})(y^{k})=\mathcal{A}(Z^{k})$,
which implies \eqref{eq:direct-y-projection}.
Using \eqref{eq:shifted-direct-multiplier} and
\eqref{eq:direct-y-projection}, the direct $S$-subproblem, up to terms independent
of $S$, becomes
\begin{equation}
    \min_{S\in\mathcal{M}}
    \left\{
        \left\langle
            D+W^{k},
            Z
        \right\rangle
        +
        \frac{\sigma_k}{2}
        \left\|
            \mathcal{P}(Z^{k})-Z
        \right\|^{2}
    \right\}.
    \label{eq:direct-ALM-subproblem}
\end{equation}
Since the ranges of $\mathcal{P}$ and $\mathcal{Q}$ are orthogonal, we have
\begin{align}
    \left\|
        \mathcal{P}(Z^{k})-Z
    \right\|^{2}
    &=
    \left\|
        \mathcal{P}(Z-Z^{k})
    \right\|^{2}
    +
    \left\|
        \mathcal{Q}(Z)
    \right\|^{2}.
    \label{eq:penalty-decomposition}
\end{align}
Therefore, for the same shifted multiplier, the direct subproblem
\eqref{eq:direct-ALM-subproblem} is exactly the reduced subproblem
\eqref{eq:reduced-ALM-subproblem} plus the additional term
\begin{equation}
    \frac{\sigma_k}{2}
    \left\|
        \mathcal{P}(Z-Z^{k})
    \right\|^{2}.
    \label{eq:additional-proximal-term}
\end{equation}
The term \eqref{eq:additional-proximal-term} does not penalize infeasibility.
Indeed, the component $\mathcal{P}(Z)$ lies in
$\operatorname{range}(\mathcal{A}^{*})$ and can be represented exactly by changing
$y$. Nevertheless, the direct formulation forces $\mathcal{P}(Z)$ to remain close
to the stale value $\mathcal{P}(Z^{k})$ inherited from the preceding $y$-update.
In contrast, in \eqref{dsdp1}, each trial value of $S$ is accompanied by its exact best response $y(S)=(\mathcal{A}\mathcal{A}^{*})^{-1}\mathcal{A}(S+C)$.
Thus, eliminating $y$ removes the block-coordinate lag between $S$ and $y$ in the augmented Lagrangian.

The additional regularization also acts in directions that are relevant to objective reduction. Since $D\in\operatorname{range}(\mathcal{A}^{*})$,
\begin{equation}
    \left\langle D,Z\right\rangle
    =
    \left\langle D,\mathcal{P}(Z)\right\rangle.
    \label{eq:objective-projector}
\end{equation}
Hence the objective depends on the $\mathcal{P}(Z)$ component, while the extra term \eqref{eq:additional-proximal-term} directly damps changes in that component.
When $\sigma_k$ is increased to improve feasibility, the direct formulation simultaneously strengthens this artificial damping. In \eqref{dsdp1}, increasing $\sigma_k$ only strengthens the penalty on the genuine feasibility residual
$\mathcal{Q}(Z)$.
This distinction can also be seen from the curvature of the two subproblems in the ambient variable $S$. Their quadratic terms have Hessians
\begin{equation}
    \nabla_S^{2}\Phi_k
    =
    \sigma_k\mathcal{Q},
    \qquad
    \nabla_S^{2}\Phi_k^{\mathrm{direct}}
    =
    \sigma_k \mathcal{I}
    =
    \sigma_k\mathcal{Q}
    +
    \sigma_k\mathcal{P}.
    \label{eq:hessian-comparison}
\end{equation}
Thus, the direct formulation adds curvature $\sigma_k\mathcal{P}$ in directions
that can already be handled exactly through $y$. After the Burer--Monteiro
factorization $S=YY^{\intercal}$, this additional curvature is inherited by the
Riemannian subproblem and can make the trust-region model unnecessarily stiff,
leading to smaller effective steps and more outer iterations in the direct ADMM scheme.

The preceding argument does not imply that \eqref{dsdp1} must be faster for every instance, since proximal regularization may occasionally improve stability.
Nevertheless, it identifies a principal mechanism behind the empirical advantage of \eqref{dsdp1}: the reformulation penalizes only the true constraint violation, updates the eliminated variable exactly for every candidate $S$, and avoids
unnecessary damping in the objective-relevant subspace. This is consistent with
the numerical results in
Section~\ref{sec:reformulation-comparison-experiments}, where the performance
gap becomes particularly pronounced as the problem size increases.

\section{Exploiting the manifold structure}\label{manifold}
Assumption~\ref{assump1} motivates a Burer--Monteiro factorization of the PSD slack variable $S$. For a factorization size $p$, define $\Psi_k(Y)\coloneqq\Phi_k(YY^{\intercal})$ for $Y\in\R^{n\times p}$.
By doing so, the convex subproblem \eqref{subp1} becomes the non-convex factorized subproblem:
\begin{equation}\label{subpY}
\min_{Y\in\cN_p}\Psi_k(Y)=L_{\sigma_{k}}(YY^{\intercal},y(YY^{\intercal}),\widetilde{X}^{k}),\tag{Sub-$k$}
\end{equation}
where 
\begin{equation}\label{maniY}
\cN_p\coloneqq\{Y\in\R^{n\times p}\mid YY^{\intercal}\in\cM\}=\{Y\in\R^{n\times p}\mid \|Y(i,:)\|=1,\, i=1,\ldots,n\},
\end{equation}
is the oblique manifold ($Y(i,:)$ stands for the $i$-th row of $Y$).
We call $p$ the \emph{factorization size}.
Thus, it remains to solve \eqref{subpY} on the manifold $\cN_p$, which could be done by off-the-shelf efficient Riemannian optimization methods. Here, we choose the Riemannian trust-region method to solve \eqref{subpY} due to its superior performance, as shown in \cite{wang2025solving}. We refer the reader to \cite{AbsBakGal2007-FoCM} for details on the Riemannian trust-region method.

We next calculate the Riemannian gradient and Riemannian Hessian of \eqref{subpY} that are necessary to perform optimization on the manifold $\cN_p$ via the Riemannian trust-region method. First of all, we note that at a point $Y\in\cN_p$, the tangent space of $\cN_p$ is given by
\begin{equation}
    T_Y\cN_p=\{U\in\R^{n\times p}\mid\diag(UY^{\intercal})=0\},
\end{equation}
and the normal space to $\cN_p$ is given by
\begin{equation}
    N_Y\cN_p=\{\Diag(u)Y\mid u\in\R^n\}.
\end{equation}

\begin{lemma}[\cite{wang2025solving}, Lemma 4.1]\label{sec4:lm1}
Let $Y$ be a point on $\cN_p$. The orthogonal projector $P_Y:\R^{n\times p}\rightarrow T_Y\cN_p$ is given by
\begin{equation}
    P_Y(U)=U-\Diag(UY^{\intercal})Y, \quad\text{ for }U\in\R^{n\times p}.
\end{equation}
\end{lemma}

\begin{proposition}\label{prop}
Consider the non-convex subproblem \eqref{subpY}. Let $S=YY^{\intercal}$ and $X=\nabla\Phi_k(S)-\Diag(\nabla\Phi_k(S)S)$. Then the \emph{Riemannian gradient} at $Y$ is given by
\begin{equation}\label{sec4:eq6}
\grad\,\Psi_k(Y)=2XY.
\end{equation}
For $U\in T_{Y}\cN_p$, let 
\begin{equation*}
 \widetilde{H}=\nabla^2\Psi_k(Y)[U]=2\nabla\Phi_k(S)U+2\sigma_k\mathcal Q(Z)Y,
\end{equation*}
where $Z\coloneqq YU^{\intercal}+UY^{\intercal}$.
Then the \emph{Riemannian Hessian} is given by
\begin{equation}\label{sec4:eq7}
\hess\,\Psi_k(Y)[U]=\widetilde{H}-\Diag(\widetilde{H}Y^{\intercal})Y-2\Diag(\nabla\Phi_k(S)S)U.
\end{equation}
\end{proposition}
\begin{proof}
By (3) of \cite{absil2013extrinsic}, $\grad\,\Psi_k(Y)=P_Y(\nabla\Psi_k(Y))=P_Y(2\nabla\Phi_k(S)Y)$. By Lemma \ref{sec4:lm1}, there exists $z\in\R^n$ such that
\begin{equation*}
\grad\,\Psi_k(Y)=2\nabla\Phi_k(S)Y-2\Diag(z)Y.
\end{equation*}
Moreover, since $\grad\,\Psi_k(Y)\in T_Y\cN_p$, we have
\begin{equation}\label{sec4:eq4}
\diag\left(\grad\,\Psi_k(Y)Y^{\intercal}\right)=2\diag\left((\nabla\Phi_k(S)-\Diag(z))S\right)=0,
\end{equation}
which gives $z=\diag(\nabla\Phi_k(S)S)$.

By (10) of \cite{absil2013extrinsic}, we have
\begin{equation}
    \hess\,\Psi_k(Y)[U]=P_Y(\nabla^2\Psi_k(Y)[U])+\mathfrak{A}_Y(U,P_Y^{\perp}(\nabla\Psi_k(Y))),
\end{equation}
where $\mathfrak{A}_Y$ is the Weingarten map at $Y$ and $P_Y^{\perp}=I-P_Y$ is the orthogonal projector at $Y$ to $N_Y\cN_p$. Let $D_Y(\cdot)[U]$ be the directional derivative at $Y$ along $U$.
Then,
\begin{align*}
\mathfrak{A}_Y(U,P_Y^{\perp}(\nabla\Psi_k(Y)))&=\mathfrak{A}_Y(U,2\Diag(z)Y)\\
&=-P_Y(D_Y(2\Diag(z)Y)[U])\\
&=-2P_Y(\Diag(z)U)-2P_Y(\Diag(D_{Y}(z)[U])Y)\\
&=-2\Diag(z)U,
\end{align*}
where we have used the fact that $P_Y(\Diag(z)U)=\Diag(z)U$ and $P_Y(\Diag(u)Y)=0$ for any $u\in\R^n$. \eqref{sec4:eq7} then follows.
\end{proof}

The optimality of $S=YY^{\intercal}$ to the convex subproblem \eqref{subp1} is determined solely by the PSDness of the certificate $X$.
\begin{proposition}\label{sec4:thm1}
Let $Y\in\cN_p$, $S=YY^{\intercal}$, and $X=\nabla\Phi_k(S)-\Diag(\nabla\Phi_k(S)S)$. Then a stationary point $Y$ of the non-convex factorized subproblem \eqref{subpY} provides a minimizer $S=YY^{\intercal}$ of the convex subproblem \eqref{subp1} if and only if $X\succeq0$.
\end{proposition}
\begin{proof}
The condition that $Y$ is a stationary point means $\grad\,\Psi_k(Y)=2XY=0$ implying $XS=0$. The conclusion then follows from Lemma \ref{lm:suboptimality}.
\end{proof}

We then obtain the following theorem as an immediate corollary of Proposition \ref{sec4:thm1}, which connects subproblem optimality to the original problem.
\begin{corollary}\label{sec4:thm3}
A stationary point $Y\in\cN_p$ of the non-convex factorized subproblem \eqref{subpY} provides a minimizer $S=YY^{\intercal}$ of \eqref{dsdp1} if and only if $\mathcal{A}^*(y(S))=S+C$ and $X=\nabla\Phi_k(S)-\Diag(\nabla\Phi_k(S)S)\succeq0$.
\end{corollary}

Because the subproblem \eqref{subpY} is nonconvex, a mechanism is needed to escape stationary points whose certificate is not PSD. The next result provides an explicit negative-curvature direction.
\begin{theorem}\label{sec5:thm1}
Suppose that $X=\nabla\Phi_k(S)-\Diag(\nabla\Phi_k(S)S)\nsucceq0$. Let $\delta\geq1$ be an integer, and let $V\in\R^{n\times\delta}$ have columns chosen from eigenvectors of $X$ associated with negative eigenvalues. Define
\[
    \widehat Y\coloneqq[Y,0_{n\times\delta}]\in\cN_{p+\delta},
    \qquad
    U\coloneqq[0_{n\times p},V].
\]
Then $U\in T_{\widehat Y}\cN_{p+\delta}$ and
\[
    \langle U,\grad\,\Psi_k(\widehat Y)\rangle=0,
    \qquad
    \langle U,\hess\,\Psi_k(\widehat Y)[U]\rangle<0.
\]
Thus, $U$ is a second-order descent direction for the rank-augmented subproblem.
\end{theorem}
\begin{proof}
The zero-column construction gives $\widehat YU^{\intercal}=0$, and hence $U\in T_{\widehat Y}\cN_{p+\delta}$. It also yields
\[
    \langle U,\grad\,\Psi_k(\widehat Y)\rangle
    =2\Tr(U^{\intercal}X\widehat Y)=0.
\]
Since $\widehat YU^{\intercal}+U\widehat Y^{\intercal}=0$, Proposition~\ref{prop} gives $\hess\,\Psi_k(\widehat Y)[U]=2XU$. Therefore,
\[
    \langle U,\hess\,\Psi_k(\widehat Y)[U]\rangle
    =2\Tr(U^{\intercal}XU)
    =2\Tr(V^{\intercal}XV)<0.
\]
\end{proof}

\section{The {\tt ManiDSDP} algorithm}\label{algorithms}
Our dual Riemannian ALM algorithm for solving \eqref{dsdp1} is presented in Algorithm~\ref{alg1}.

\begin{algorithm}
\renewcommand{\algorithmicrequire}{\textbf{Input:}}
\renewcommand{\algorithmicensure}{\textbf{Output:}}
\caption{{\tt ManiDSDP}}\label{alg1}
	\begin{algorithmic}[1]
   	\REQUIRE $\cA,b,C$; $\sigma_0>0$; $0<\sigma_{\min}<\sigma_{\max}$; $p_0\geq1$; $\gamma>1$; $0<\kappa_1\leq\kappa_2$
		\STATE $k \gets 0$, $p\gets p_0$, $\widetilde{X}^0\gets 0_{n\times n}$, $Y^0 \gets 0_{n\times p}$, $U \gets 0_{n\times p}$;
		\WHILE{The stopping criterion is not fulfilled}
		\STATE Approximately solve \eqref{subpY} by a Riemannian trust-region method, warm-started at $Y^k$ and supplied with the carried descent direction $U$, to obtain $Y^{k+1}\in\cN_p$;
		\STATE $S^{k+1} \gets Y^{k+1}(Y^{k+1})^{\intercal}$; 
        \STATE $y^{k+1} \gets (\cA\cA^*)^{-1}\cA(S^{k+1}+C)$;
        \STATE $\widetilde{X}^{k+1}\gets \widetilde{X}^k-\sigma_k\left(\cA^*\left(y^{k+1}\right)-S^{k+1}-C\right)$;
        \STATE $z^{k+1}\gets\diag\left(\left(\widetilde{X}^{k+1}+D\right)S^{k+1}\right)$;
        \STATE $X^{k+1}\gets\widetilde{X}^{k+1}+D-\Diag\left(z^{k+1}\right)$;
		\STATE If $\lambda_{\min}(X^{k+1})<0$, form a second-order descent direction $U$ as in Theorem~\ref{sec5:thm1} and set $U\gets0$ otherwise;
		\STATE Update $p$ following the strategy described in \cite[Section~5.1]{wang2025solving}, padding $Y^{k+1}$ and $U$ with zero columns as needed;
        \STATE Update $\sigma_k$ according to the rules \eqref{sec5:eq1}; 
        \STATE $k \gets k+1$;
		\ENDWHILE
		\RETURN $(S^{k},y^{k},X^{k},z^{k})$;
	\end{algorithmic}
\end{algorithm}

We make a few remarks about Algorithm \ref{alg1}:
\begin{itemize}
    \item Most matrix operations can be expressed in terms of $Y^{k+1}$, avoiding explicit formation of the dense matrix $S^{k+1}=Y^{k+1}(Y^{k+1})^{\intercal}$ whenever possible.
    \item The algorithm does not require forming a dense inverse at every iteration. One factors $\mathcal A\mathcal A^*$ once and reuses the factorization to apply $(\mathcal A\mathcal A^*)^{-1}$ when evaluating $y(S)$ and the orthogonal projector. For SDPs arising from moment--SOS relaxations, $\mathcal A$ has orthogonal or partially orthogonal coefficient structure, so $\mathcal A\mathcal A^*$ is diagonal or admits a special structure that can substantially simplify the computation of $(\mathcal A\mathcal A^*)^{-1}$ (\cite[Theorem A11]{yang2020one}).
    \item The initial factorization size $p_0$ can be set to $2$. However, we have found that choosing $p_0\sim O(\lceil\log(m)\rceil)$ would sometimes bring some speed-up; see \cite{so2008unified}.
    \item In Step 9, a partial eigenvalue decomposition could be performed when the full eigenvalue decomposition is too expensive.
    \item The direction $U$ computed in Step~9 is passed to the next execution of Step~3. Therefore, it is a second-order descent direction for the subproblem just completed and a \emph{one-step-delayed} curvature direction for the next subproblem. Section~\ref{subsec:one-step-delay} quantifies the resulting perturbation.
    \item To minimize the cost of solving the ALM subproblem \eqref{subpY}, we follow the strategy described in \cite[Section~5.1]{wang2025solving} to dynamically adjust the factorization size $p$; 
    \item The following self-adaptive strategy is adopted to update the penalty factor $\sigma$:
    \begin{equation}\label{sec5:eq1}
    \sigma_{k+1}=\begin{cases}
    \max\,\{\sigma_k/\gamma,\sigma_{\min}\},&\text{ if }\|R^{k+1}\|/(1 +\left\lVert b \right\rVert)<\kappa_1\|\grad\,\Psi_k(Y^{k+1})\|,\\
    \min\,\{\gamma\sigma_k,\sigma_{\max}\},&\text{ if }\|R^{k+1}\|/(1 +\left\lVert b \right\rVert)>\kappa_2\|\grad\,\Psi_k(Y^{k+1})\|,\\
    \sigma_k,&\text{ otherwise,}\\
    \end{cases}
     \end{equation}
     where $R^{k+1}=\mathcal{A}^*(y^{k+1})-S^{k+1}-C$.
     \item Algorithm \ref{alg1} can be extended to handle SDPs with multi-blocks which often arise as relaxations of sparse polynomial optimization problems. 
\end{itemize}

\subsection{The one-step-delayed design}\label{subsec:one-step-delay}
For $S\in\mathcal M$, define the linear certificate map
\begin{equation}\label{eq:general-certificate-map}
    \mathscr C_S(G)
    \coloneqq
    G-\Diag\bigl(\diag(GS)\bigr),
    \qquad G\in\mathbb S_n.
\end{equation}
At the end of outer iteration $k$, one has
\begin{equation}\label{eq:completed-certificate}
    X^{k+1}
    =\mathscr C_{S^{k+1}}(\widetilde X^{k+1}+D).
\end{equation}
This is the certificate of the completed subproblem $\Psi_k$. At the same matrix $S^{k+1}$, the gradient of the next subproblem is
\[
\nabla\Phi_{k+1}(S^{k+1})
=\widetilde X^{k+1}+D-\sigma_{k+1}R^{k+1}.
\]
Consequently, the certificate relevant to the starting point of $\Psi_{k+1}$ is
\begin{equation}\label{eq:delayed-certificate}
    \widehat X^{k+1}
    \coloneqq
    \mathscr C_{S^{k+1}}\bigl(\nabla\Phi_{k+1}(S^{k+1})\bigr)
    =X^{k+1}-\sigma_{k+1}\mathscr C_{S^{k+1}}(R^{k+1}).
\end{equation}
Thus, the one-step delay changes the certificate by a term proportional to the current feasibility residual and to the newly selected penalty parameter.

\begin{proposition}\label{prop:delayed-curvature}
Suppose that $\nu_{k+1}\coloneqq-\lambda_{\min}(X^{k+1})>0$, and let $v^{k+1}$ be a unit eigenvector satisfying $X^{k+1}v^{k+1}=-\nu_{k+1}v^{k+1}$.
After appending one zero column, set
\[
    \widehat Y^{k+1}\coloneqq[Y^{k+1},0],
    \qquad
    U^{k+1}\coloneqq[0,v^{k+1}].
\]
Then $U^{k+1}\in T_{\widehat Y^{k+1}}\mathcal N_{p+1}$ and
\begin{equation}\label{eq:delayed-first-order-orthogonality}
    \big\langle U^{k+1},\grad\,\Psi_{k+1}(\widehat Y^{k+1})\big\rangle=0.
\end{equation}
Moreover,
\begin{align}
\big\langle U^{k+1},
\hess\,\Psi_{k+1}(\widehat Y^{k+1})[U^{k+1}]\big\rangle
&=-2\nu_{k+1}
-2\sigma_{k+1}(v^{k+1})^{\intercal}
\mathscr C_{S^{k+1}}(R^{k+1})v^{k+1} \notag\\
&\leq
-2\nu_{k+1}
+2(1+\sqrt n)\sigma_{k+1}\|R^{k+1}\|.                 \label{eq:delayed-curvature-bound}
\end{align}
In particular, the carried direction remains a strict negative-curvature direction for $\Psi_{k+1}$ whenever
\begin{equation}\label{eq:curvature-persistence-test}
    \nu_{k+1}>
    (1+\sqrt n)\sigma_{k+1}\|R^{k+1}\|.
\end{equation}
\end{proposition}

\begin{proof}
The zero-column construction gives
$\widehat Y^{k+1}(U^{k+1})^{\intercal}=0$, and hence
$U^{k+1}\in T_{\widehat Y^{k+1}}\mathcal N_{p+1}$.
It also makes the nonzero columns of $U^{k+1}$ disjoint from those of
$\widehat Y^{k+1}$, which proves
\eqref{eq:delayed-first-order-orthogonality}. Applying the calculation in the
proof of Theorem~\ref{sec5:thm1} to the objective $\Psi_{k+1}$ gives
\[
\big\langle U^{k+1},
\hess\,\Psi_{k+1}(\widehat Y^{k+1})[U^{k+1}]\big\rangle
=2(v^{k+1})^{\intercal}\widehat X^{k+1}v^{k+1}.
\]
Substitution of \eqref{eq:delayed-certificate} yields the equality in
\eqref{eq:delayed-curvature-bound}.
It remains to bound the perturbation. For $S\in\mathcal M$ and the $i$-th unit
vector $e_i$,
\[
\|Se_i\|^2=e_i^{\intercal}S^2e_i
\leq\|S\|_2e_i^{\intercal}Se_i
=\|S\|_2
\leq\Tr(S)=n.
\]
Therefore,
\[
|(GS)_{ii}|
\leq\|G\|_2\|Se_i\|
\leq\sqrt n\,\|G\|_2,
\]
and hence
\begin{equation}\label{eq:certificate-map-bound}
    \|\mathscr C_S(G)\|_2
    \leq(1+\sqrt n)\|G\|_2
    \leq(1+\sqrt n)\|G\|.
\end{equation}
The inequality in \eqref{eq:delayed-curvature-bound} and the sufficient
condition \eqref{eq:curvature-persistence-test} follow.
\end{proof}

The one-step-delayed design of {\tt ManiDSDP} avoids repeatedly reopening a subproblem that is about to be
replaced by a new augmented Lagrangian model. It also amortizes rank expansion and
eigenvalue computations over the outer ALM iterations, and no vector transport is needed
because the old factor is the starting factor of the next subproblem. Its limitation is
that the residual term in \eqref{eq:delayed-certificate}, especially after an increase
of $\sigma_{k+1}$, can make the direction stale or reverse its curvature sign. Thus the
direction $U$ in Algorithm~\ref{alg1} should be understood as a curvature-informed warm start. Merely passing it to the next subproblem does not by itself imply global
convergence; a quantitative condition on the spectral defect produced by the next iteration is required.
We next formulate such a condition and establish global convergence.

\subsection{Global convergence}
For $k\geq1$, define the negative-eigenvalue defect and its penalty-weighted
version by
\begin{equation}\label{eq:delayed-defects}
    \nu_k\coloneqq\max\{0,-\lambda_{\min}(X^k)\},
    \qquad
    d_k\coloneqq2n\sigma_{k-1}\nu_k.
\end{equation}
The defect $d_k$ is detected after completing $\Psi_{k-1}$, and the direction
constructed from $X^k$ is used when solving $\Psi_k$, which produces $d_{k+1}$.

\begin{assumption}\label{assump:delayed-forcing}
There exist constants $0\leq\rho<1$ and $0\leq\theta<1-\rho$, together with a
nonnegative sequence $\{\varepsilon_k\}_{k\geq1}$ satisfying
$\sum_{k=1}^{\infty}\varepsilon_k<\infty$, such that
\begin{equation}\label{eq:delayed-forcing}
    d_{k+1}
    \leq
    \rho d_k
    +\theta\sigma_k^2\|R^{k+1}\|^2
    +\varepsilon_k,
    \qquad k\geq1.
\end{equation}
\end{assumption}

Assumption~\ref{assump:delayed-forcing} is a one-step-delayed inexactness
condition rather than an additional algorithmic step. It allows a fraction $\rho$ of
the incoming spectral defect to remain, allows a residual-controlled error that can
be absorbed by the augmented Lagrangian decrease, and allows summable errors from solving the inner subproblem. We now establish the global convergence of {\tt ManiDSDP} under Assumption~\ref{assump:delayed-forcing}. Let us begin by proving a useful lemma.

\begin{lemma}\label{sec5:lm1}
Let $(S^{\star},y^{\star},X^{\star},z^{\star})$ be a KKT point of
\eqref{dsdp1}, set
$\widetilde X^{\star}=X^{\star}-D+\mathcal B^*(z^{\star})$, and let
$R^{k+1}=\mathcal A^*(y^{k+1})-S^{k+1}-C$. Then
\begin{equation}\label{sec5:eq5}
    \langle\widetilde X^k-\widetilde X^{\star},R^{k+1}\rangle
    \geq
    \sigma_k\|R^{k+1}\|^2-n\nu_{k+1}.
\end{equation}
\end{lemma}

\begin{proof}
Since $X^{k+1}\succeq-\nu_{k+1}I$, $S^{\star}\succeq0$, and
$\Tr(S^{\star})=n$, Lemma~\ref{sec3:lm1} gives
\begin{equation}\label{sec5:eq2}
    \langle X^{k+1},S^{\star}-S^{k+1}\rangle
    =\langle X^{k+1},S^{\star}\rangle
    \geq-n\nu_{k+1}.
\end{equation}
Moreover,
\[
S^{\star}-S^{k+1}
=\mathcal A^*(y^{\star}-y^{k+1})+R^{k+1}.
\]
Using $\mathcal A(\widetilde X^{k+1})=0$ and
$\widetilde X^{k+1}=\widetilde X^k-\sigma_kR^{k+1}$, we obtain
\begin{equation}\label{eq:delayed-lemma-xtilde}
\langle\widetilde X^{k+1},S^{\star}-S^{k+1}\rangle
=\langle\widetilde X^k,R^{k+1}\rangle
-\sigma_k\|R^{k+1}\|^2.
\end{equation}
Moreover, using
$\mathcal A(D)=b=\mathcal A(X^{\star}+\mathcal B^*(z^{\star}))$,
$\langle X^{\star},S^{\star}\rangle=0$,
$\langle X^{\star},S^{k+1}\rangle\geq0$, and
$\langle\mathcal B^*(z^{\star}),S^{\star}-S^{k+1}\rangle=0$, we obtain
\begin{align}
\langle D,S^{\star}-S^{k+1}\rangle
&=\langle D,\mathcal A^*(y^{\star}-y^{k+1})+R^{k+1}\rangle\notag\\
&=\langle \mathcal A(D),y^{\star}-y^{k+1}\rangle
  +\langle D,R^{k+1}\rangle\notag\\
&=\langle X^{\star}+\mathcal B^*(z^{\star}),
       \mathcal A^*(y^{\star})-\mathcal A^*(y^{k+1})\rangle
  +\langle D,R^{k+1}\rangle\notag\\
&=\langle X^{\star}+\mathcal B^*(z^{\star}),
       S^{\star}-S^{k+1}-R^{k+1}\rangle
  +\langle D,R^{k+1}\rangle\notag\\
&\leq\langle D-X^{\star}-\mathcal B^*(z^{\star}),R^{k+1}\rangle\notag\\
&=-\langle\widetilde X^{\star},R^{k+1}\rangle.
\label{eq:delayed-lemma-D}
\end{align}
Finally, both $S^{\star}$ and $S^{k+1}$ have unit diagonal, so
\[
\langle-\mathcal B^*(z^{k+1}),S^{\star}-S^{k+1}\rangle=0.
\]
Substituting
$X^{k+1}=\widetilde X^{k+1}+D-\mathcal B^*(z^{k+1})$ into
\eqref{sec5:eq2} and using
\eqref{eq:delayed-lemma-xtilde}--\eqref{eq:delayed-lemma-D} proves
\eqref{sec5:eq5}.
\end{proof}

\begin{theorem}[Global convergence]
\label{thm:multiplier-full-sequence-convergence}
Let $\{(S^k,y^k,X^k,z^k,\widetilde X^k)\}_{k\geq1}$ be generated by
Algorithm~\ref{alg1}. Assume that the KKT set of \eqref{dsdp} is nonempty,
$\sigma_{\min}\leq\sigma_k\leq\sigma_{\max}$ for all $k$, and
Assumption~\ref{assump:delayed-forcing} holds. Then:
\begin{enumerate}[(i)]
    \item $\lim_{k\to\infty}R^k=0$ and $\sum_{k=1}^{\infty}\nu_k<\infty$;
    \item The sequences
    $\{\widetilde X^k\}$, $\{S^k\}$, $\{y^k\}$, $\{X^k\}$, and
    $\{z^k\}$ are bounded;
    \item Every cluster point of
    $\{(S^k,y^k,X^k,z^k)\}$ is a KKT point of \eqref{dsdp};
    \item There exists a KKT point
    $(\overline S,\overline y,\overline X,\overline z)$ such that, with
    $\overline{\widetilde X}
    \coloneqq\overline X-D+\Diag(\overline z)$,
    \[
        \widetilde X^k\longrightarrow\overline{\widetilde X},
        \qquad
        \operatorname{dist}\bigl((S^k,y^k,X^k,z^k),\mathcal K\bigr)
        \longrightarrow0,
    \]
    where $\mathcal K$ is the KKT set of \eqref{dsdp}.
\end{enumerate}
\end{theorem}

\begin{proof}
Fix an arbitrary KKT point
$(S^{\star},y^{\star},X^{\star},z^{\star})$ and set
$\widetilde X^{\star}=X^{\star}-D+\mathcal B^*(z^{\star})$. Define $a_k\coloneqq\|\widetilde X^k-\widetilde X^{\star}\|^2$.
By Lemma~\ref{sec5:lm1} and the multiplier update,
\begin{align}
    a_{k+1}\leq
    a_k-\sigma_k^2\|R^{k+1}\|^2
    +2n\sigma_k\nu_{k+1}=a_k-\sigma_k^2\|R^{k+1}\|^2+d_{k+1}.
    \label{eq:delayed-basic-fejer}
\end{align}
Set
\begin{equation}\label{eq:delayed-lyapunov-constants}
    c\coloneqq\frac{\rho}{1-\rho},
    \qquad
    \alpha\coloneqq\frac{1-\rho-\theta}{1-\rho}>0.
\end{equation}
Since $1+c=(1-\rho)^{-1}$ and $(1+c)\rho=c$, combining
\eqref{eq:delayed-basic-fejer} with \eqref{eq:delayed-forcing} gives
\begin{align}
    a_{k+1}+cd_{k+1}
    &\leq
    a_k+cd_k
    -\alpha\sigma_k^2\|R^{k+1}\|^2
    +\frac{\varepsilon_k}{1-\rho}.
    \label{eq:delayed-lyapunov-step}
\end{align}
Define
\begin{equation}\label{eq:delayed-lyapunov}
    \mathcal V_k
    \coloneqq
    a_k+cd_k+\frac{1}{1-\rho}
    \sum_{j=k}^{\infty}\varepsilon_j.
\end{equation}
Then
\begin{equation}\label{eq:delayed-lyapunov-descent}
    \mathcal V_{k+1}
    \leq
    \mathcal V_k
    -\alpha\sigma_k^2\|R^{k+1}\|^2.
\end{equation}
Thus, $\{\mathcal V_k\}$ is nonincreasing and bounded below. In particular,
$\{\widetilde X^k\}$ is bounded, and summing
\eqref{eq:delayed-lyapunov-descent} yields
\begin{equation}\label{eq:residual-square-summability}
    \sum_{k=1}^{\infty}\|R^{k+1}\|^2<\infty.
\end{equation}
Hence $\lim_{k\to\infty}R^k=0$.

Summing \eqref{eq:delayed-forcing} from $k=1$ to $N$ gives
\begin{align*}
    (1-\rho)\sum_{k=2}^{N}d_k+d_{N+1}
    \leq
    \rho d_1
    +\theta\sum_{k=1}^{N}\sigma_k^2\|R^{k+1}\|^2
    +\sum_{k=1}^{N}\varepsilon_k.
\end{align*}
The right-hand side is uniformly bounded by
\eqref{eq:residual-square-summability} and the summability of
$\{\varepsilon_k\}$. Consequently,
$\sum_{k=1}^{\infty}d_k<\infty$. Since
$d_k\geq2n\sigma_{\min}\nu_k$, we obtain
$\sum_{k=1}^{\infty}\nu_k<\infty$, and in particular $\lim_{k\to\infty}\nu_k=0$.
This proves (i).

For every $k$, $S^k\in\mathcal M$, and therefore
$\Tr(S^k)=n$ and $\|S^k\|\leq n$. Hence $\{S^k\}$ is bounded. The closed-form
update for $y^k$ implies that $\{y^k\}$ is bounded. Furthermore,
\[
    z^k=\diag\bigl((\widetilde X^k+D)S^k\bigr),
    \qquad
    X^k=\widetilde X^k+D-\Diag(z^k).
\]
So $\{z^k\}$ and $\{X^k\}$ are bounded as well. This proves (ii).

Let
\[
(S^{k_j},y^{k_j},X^{k_j},z^{k_j})
\longrightarrow
(\widehat S,\widehat y,\widehat X,\widehat z)
\]
be an arbitrary convergent subsequence. Since $R^{k_j}\to0$,
$\mathcal A^*(\widehat y)-\widehat S-C=0$. Closedness of $\mathcal M$ gives
$\widehat S\succeq0$ and $\diag(\widehat S)=\mathbf1$. Since
$X^k\succeq-\nu_kI$ and $\nu_k\to0$, we have $\widehat X\succeq0$.
Lemma~\ref{sec3:lm1} gives
$\langle\widehat X,\widehat S\rangle=0$, and hence
$\widehat X\widehat S=0$. Finally,
\[
    \mathcal A(X^k+\mathcal B^*(z^k))
    =\mathcal A(\widetilde X^k+D)=b
\]
for every $k$. Passing to the limit proves dual affine feasibility. Thus the
cluster point is a KKT point, proving (iii).

To prove (iv), choose a cluster point
$(\overline S,\overline y,\overline X,\overline z)$ and set
$\overline{\widetilde X}=\overline X-D+\Diag(\overline z)$. Along the
corresponding subsequence,
$\widetilde X^{k_j}\to\overline{\widetilde X}$. The Lyapunov argument above
applies with this KKT point in place of the initially chosen one. Since
$d_k\to0$ and the tail sum in \eqref{eq:delayed-lyapunov} tends to zero, the
scalar sequence
$\|\widetilde X^k-\overline{\widetilde X}\|^2$ has a finite limit. Its value
along the subsequence is zero, and therefore
$\widetilde X^k\to\overline{\widetilde X}$. The convergence to the KKT set
then follows by the usual compactness contradiction: otherwise a subsequence
would remain a fixed positive distance from $\mathcal K$, while a further
convergent subsequence would have a KKT limit by (iii).
\end{proof}

\begin{remark}
In the implementation, if the condition \eqref{eq:delayed-forcing} is not satisfied, it can be enforced by repeatedly computing a second-order descent direction, increasing the factorization size $p$, and solving the same Riemannian subproblem.
\end{remark}

\begin{corollary}\label{cor:spectral-tolerance}
Suppose that the assumptions of Theorem~\ref{thm:multiplier-full-sequence-convergence}
other than Assumption~\ref{assump:delayed-forcing} hold, and that there exists a
nonnegative sequence $\{\tau_k\}_{k\geq0}$ satisfying
\begin{equation}\label{stopcri}
    \sum_{k=0}^{\infty}\tau_k<\infty,
    \qquad
    \lambda_{\min}(X^{k+1})\geq-\tau_k
    \quad\text{for all }k\geq0.
\end{equation}
Then all conclusions of Theorem~\ref{thm:multiplier-full-sequence-convergence}
hold.
\end{corollary}

\begin{proof}
Condition \eqref{stopcri} gives $\nu_{k+1}\leq\tau_k$, and hence
\[
    d_{k+1}=2n\sigma_k\nu_{k+1}
    \leq2n\sigma_{\max}\tau_k.
\]
Thus \eqref{eq:delayed-forcing} holds with $\rho=\theta=0$ and
$\varepsilon_k=2n\sigma_{\max}\tau_k$, which is summable.
\end{proof}


\begin{remark}
Theorem~\ref{thm:multiplier-full-sequence-convergence} establishes genuine
full-sequence convergence only for $\widetilde X^k$. Without an additional
uniqueness or isolation assumption, the whole sequence
$\{S^k,y^k,X^k,z^k\}$ need not converge to a single KKT point; it is only
guaranteed to approach the KKT set.
\end{remark}

\subsection{Finite rank identification}
We next establish finite rank identification for Algorithm~\ref{alg1}. 
For a symmetric matrix $G\in\mathbb S_n$, write its eigenvalues in increasing
and decreasing order as
\[
 \lambda_1^{\uparrow}(G)\leq\cdots\leq\lambda_n^{\uparrow}(G),
 \qquad
 \lambda_1^{\downarrow}(G)\geq\cdots\geq\lambda_n^{\downarrow}(G),
\]
respectively.

\begin{theorem}[Finite rank identification]\label{thm:rank-identification}
Let $\{(Y^k,S^k,X^k)\}_{k\geq1}$ be generated by
Algorithm~\ref{alg1}, where $Y^k\in\mathbb R^{n\times p_k}$.
Suppose that the following conditions hold:
\begin{enumerate}[(i)]
    \item The whole primal-dual sequence converges to a KKT pair,
    \[
        S^k\longrightarrow S^\star,
        \qquad
        X^k\longrightarrow X^\star;
    \]

    \item The limiting pair is strictly complementary, namely,
    \[
        S^\star\succeq0,
        \qquad
        X^\star\succeq0,
        \qquad
        S^\star X^\star=0,
        \qquad
        \rank(S^\star)+\rank(X^\star)=n;
    \]

    \item There exists $K_0\in\mathbb N$ such that the exact complementarity relation $X^kY^k=0$ holds for every $k\geq K_0$.
\end{enumerate}
Then there exists $N\in\mathbb N$ such that
\begin{equation}\label{eq:finite-rank-identification}
    \rank(S^k)=\rank(S^\star)
    \qquad\text{for every }k\geq N.
\end{equation}
\end{theorem}

\begin{proof}
Let $r\coloneqq\rank(S^\star)$.
If $r=n$, then $S^\star\succ0$.  Since $S^k\to S^\star$, Weyl's
eigenvalue inequality gives
\[
    \lambda_n^{\downarrow}(S^k)
    \longrightarrow
    \lambda_n^{\downarrow}(S^\star)>0.
\]
Hence $S^k\succ0$ and $\rank(S^k)=n=r$ for all sufficiently large $k$.

Assume henceforth that $r<n$.  Strict complementarity yields
$\rank(X^\star)=n-r$.
Since $X^\star\succeq0$, it has exactly $r$ zero eigenvalues, and therefore $\alpha\coloneqq\lambda_{r+1}^{\uparrow}(X^\star)>0$.
By $X^k\to X^\star$ and Weyl's eigenvalue inequality, there exists $K_1\in\mathbb N$ such that
\begin{equation}\label{eq:dual-gap-persistence}
    \lambda_{r+1}^{\uparrow}(X^k)
    \geq \frac{\alpha}{2}>0
    \qquad\text{for every }k\geq K_1.
\end{equation}
Consequently, $X^k$ has at most $r$ nonpositive eigenvalues and, in particular,
\begin{equation}\label{eq:nullity-upper-bound}
    \dim\ker(X^k)\leq r
    \qquad\text{for every }k\geq K_1.
\end{equation}
For $k\geq\max\{K_0,K_1\}$, $X^kY^k=0$ gives
$\operatorname{range}(Y^k)\subseteq\ker(X^k)$.
Using $S^k=Y^k(Y^k)^{\intercal}$, we obtain
$\rank(S^k)=\rank(Y^k)\leq\dim\ker(X^k)\leq r$.
It remains to prove the reverse inequality. Since $S^\star\succeq0$ has rank $r$,
$\beta\coloneqq\lambda_r^{\downarrow}(S^\star)>0$.
Again, $S^k\to S^\star$ and Weyl's eigenvalue inequality imply that there exists $K_2\in\mathbb N$ such that
\begin{equation}\label{eq:primal-gap-persistence}
    \lambda_r^{\downarrow}(S^k)
    \geq\frac{\beta}{2}>0
    \qquad\text{for every }k\geq K_2.
\end{equation}
Thus, $\rank(S^k)\geq r$ for every $k\geq K_2$, 
and we conclude that
$\rank(S^k)=r=\rank(S^\star)$
for every $k\geq N\coloneqq\max\{K_0,K_1,K_2\}$.
\end{proof}

\begin{remark}\label{rem:exact-tail-criticality}
The condition $X^kY^k=0$ follows if the Riemannian
subproblem is solved exactly to first-order criticality for all sufficiently large iterations. Indeed, at iteration $k-1$,
\begin{equation*}
\nabla\Phi_{k-1}(S^k)=\widetilde X^{k-1}+D
      -\sigma_{k-1}\bigl(\mathcal A^*(y^k)-S^k-C\bigr)=\widetilde X^k.
\end{equation*}
Hence the definition of $X^k$ in Algorithm~\ref{alg1} and
Proposition~\ref{prop} give $\grad\,\Psi_{k-1}(Y^k)=2X^kY^k$.
Thus $\grad\,\Psi_{k-1}(Y^k)=0$ implies $X^kY^k=0$.
\end{remark}

\begin{remark}
The condition $X^kY^k=0$ also follows if
$X^k\succeq0$ for all sufficiently large $k$. Indeed,
Lemma~\ref{sec3:lm1} gives $\langle X^k,S^k\rangle=0$.  Since both matrices are PSD, this identity implies $X^kS^k=0$. As
$S^k=Y^k(Y^k)^{\intercal}$, one has
\[
    0=\Tr(X^kS^k)
     =\Tr\bigl((Y^k)^{\intercal}X^kY^k\bigr)
     =\bigl\|(X^k)^{1/2}Y^k\bigr\|^2,
\]
and therefore $X^kY^k=0$.
\end{remark}

\begin{remark}
Neither Assumption~\ref{assump:delayed-forcing} nor the special condition
\eqref{stopcri} implies the exact tail relation $X^kY^k=0$. Both allow small
positive eigenvalues of $S^k$ to persist at every finite iteration while
converging to zero. Thus, without exact tail criticality, eventual dual
feasibility, or an explicit hard rank-truncation step, one can generally
guarantee only asymptotic or numerical rank identification, not the exact
equality \eqref{eq:finite-rank-identification} after finitely many iterations.
\end{remark}

\section{Local one-step locking under a rank-one solution}\label{finite-locking}
In this section, we give a local explanation of the residue diving phenomenon for random BQPs observed in Section~\ref{rdp}. The key observation is that a strictly complementary rank-one solution is a sharp minimizer of all sufficiently close ALM $S$-subproblems.  Consequently, once the multiplier enters a
certain neighborhood, an exact $S$-update identifies the rank-one solution and remains there permanently.

Recall the projectors $\mathcal P$ and $\mathcal Q$ from
\eqref{eq:projection-operators}.  In view of Lemma~\ref{lm1}, the multiplier
invariant $\mathcal A(\widetilde X^k)=0$ holds for every $k\geq 1$. Thus, up to
an additive constant independent of $S$, the ALM subproblem is
\begin{equation}\label{eq:locking-subproblem}
    \min_{S\in\mathcal M}\ \Phi_{\widetilde X,\sigma}(S)
    \coloneqq
    \langle D+\widetilde X,S+C\rangle
    +\frac{\sigma}{2}\|\mathcal Q(S+C)\|^2,
\end{equation}
where $\mathcal A(\widetilde X)=0$ and $\sigma>0$. Let $s\in\{\pm1\}^n$ and define
\begin{equation}\label{eq:certificate-map}
    S^\star\coloneqq ss^{\intercal},
    \qquad
    P_s\coloneqq I-\frac{1}{n}ss^{\intercal},
    \qquad
    \mathscr C_s(G)
    \coloneqq
    G-\Diag\big((Gs)\circ s\big).
\end{equation}
Here $P_s$ is the orthogonal projector onto $s^\perp$, and $\circ$ denotes the Hadamard product. Notice
that, for every $G\in\mathbb S_n$,
\begin{equation}\label{eq:certificate-annihilates-s}
    \mathscr C_s(G)s=0,
\end{equation}
because
$\Diag((Gs)\circ s)s=Gs$.

\begin{theorem}[Local one-step locking]\label{thm:finite-locking}
Suppose that $S^\star=ss^{\intercal}$ is feasible for \eqref{dsdp1}.
Let $\widetilde X^\star\in\ker(\mathcal A)$ and define $X^\star
\coloneqq\mathscr C_s(D+\widetilde X^\star)$.
Assume the strict-complementarity condition $X^\star\succeq\delta P_s$ for some $\delta>0$. For any $\widetilde X\in\ker(\mathcal A)$ satisfying
\begin{equation}\label{eq:locking-neighborhood}
    \varepsilon(\widetilde X)
    \coloneqq
    \big\|\mathscr C_s(\widetilde X-\widetilde X^\star)\big\|_2
    <\delta,
\end{equation}
set $\gamma(\widetilde X)\coloneqq\delta-\varepsilon(\widetilde X)>0$.  Then,
for every $\sigma>0$, $S^\star$ is the unique global minimizer of
\eqref{eq:locking-subproblem}. Consequently, suppose that the $S$-subproblem in Algorithm~\ref{alg1} is solved exactly and that the generated multipliers satisfy $\widetilde X^k\to\widetilde X^\star$. Then there is an index $K<\infty$ such that
\begin{equation}\label{eq:finite-locking-conclusion}
    S^{k+1}=S^\star
    \qquad\text{for every }k\geq K.
\end{equation}
After the first such update, the feasibility residual is zero, the multiplier does not change, and the algorithm stays locked at a KKT point of \eqref{dsdp}.
\end{theorem}

\begin{proof}
For a current multiplier $\widetilde X$, let
\begin{equation}\label{eq:current-locking-certificate}
    G\coloneqq \widetilde X+D,
    \qquad
    z\coloneqq(Gs)\circ s,
    \qquad
    X\coloneqq\mathscr C_s(G)=G-\Diag(z).
\end{equation}
By \eqref{eq:certificate-annihilates-s}, both $X^\star$ and $X$ annihilate
$s$.  Moreover, the linearity of $\mathscr C_s$ gives
\begin{equation}\label{eq:certificate-perturbation}
    X-X^\star
    =\mathscr C_s(\widetilde X-\widetilde X^\star).
\end{equation}
Take an arbitrary $u\in\mathbb R^n$ and write its orthogonal decomposition as
$u=\alpha s+v$, where $v\perp s$.  Since $Xs=X^\star s=0$, it follows from \eqref{eq:certificate-perturbation}
that
\begin{align*}
    u^{\intercal}Xu
    &=v^{\intercal}Xv\\
    &=v^{\intercal}X^\star v
      +v^{\intercal}(X-X^\star)v\\
    &\geq
      \delta\|v\|^2-\|X-X^\star\|_2\|v\|^2\\
    &=\gamma(\widetilde X)\|v\|^2
     =\gamma(\widetilde X)u^{\intercal}P_su.
\end{align*}
As this holds for every $u$, we have
\begin{equation}\label{eq:current-strict-complementarity}
    X\succeq\gamma(\widetilde X)P_s\succeq0,
    \qquad
    \ker(X)=\operatorname{span}\{s\}.
\end{equation}
The kernel assertion follows because $Xs=0$ and, if $Xu=0$, then
\eqref{eq:current-strict-complementarity} yields $P_su=0$.

By the feasibility 
of $S^\star$, for any $S\in\mathcal M$ we have the exact
identity
\begin{align}
    \Phi_{\widetilde X,\sigma}(S)
    -\Phi_{\widetilde X,\sigma}(S^\star)
    ={}&\langle G,S-S^\star\rangle
    +\frac{\sigma}{2}\|\mathcal Q(S+C)\|^2.       \label{eq:objective-difference}
\end{align}
The diagonal constraints on $S$ and $S^\star$ imply
\begin{equation*}
    \langle\Diag(z),S-S^\star\rangle
    =z^{\intercal}\big(\diag(S)-\diag(S^\star)\big)=0.
\end{equation*}
Furthermore, $Xs=0$ gives
$\langle X,S^\star\rangle=s^{\intercal}Xs=0$.  Since
$G=X+\Diag(z)$, \eqref{eq:objective-difference} therefore becomes
\begin{equation}\label{eq:exact-gap-identity}
    \Phi_{\widetilde X,\sigma}(S)
    -\Phi_{\widetilde X,\sigma}(S^\star)
    =\langle X,S\rangle
    +\frac{\sigma}{2}\|\mathcal Q(S+C)\|^2.
\end{equation}
Because $S\succeq0$ and $X\succeq0$, we deduce from \eqref{eq:exact-gap-identity} that $S^\star$ is a global minimizer of \eqref{eq:locking-subproblem}.
It remains to show the uniqueness.
Let $\widehat S\coloneqq \Diag(s)S\Diag(s)$.  Then
$\widehat S\succeq0$, $\diag(\widehat S)=\mathbf 1$, and
$|\widehat S_{ij}|\leq1$ for all $i,j$.  In particular,
$1-\widehat S_{ij}\geq0$. Since multiplication by $\Diag(s)$ preserves the
Frobenius norm, we have
\begin{align}
    \|S-S^\star\|
    &=\|\widehat S-\mathbf 1\mathbf 1^{\intercal}\|\notag\leq\sum_{i,j=1}^n(1-\widehat S_{ij})\notag\\
    &=n^2-\mathbf 1^{\intercal}\widehat S\mathbf 1
      =n^2-s^{\intercal}Ss
      =n\langle P_s,S\rangle.            \label{eq:linear-distance-bound}
\end{align}
Because $S\succeq0$ and $X-\gamma(\widetilde X)P_s\succeq0$, we have $\langle X,S\rangle\geq\gamma(\widetilde X)\langle P_s,S\rangle$.
If $S\neq S^\star$, then
\eqref{eq:linear-distance-bound} shows that the
objective gap is strictly positive. Thus, $S^\star$ is the unique global
minimizer of \eqref{eq:locking-subproblem}. 

We finally prove local one-step locking.  Since $\widetilde X^k\to\widetilde X^\star$
and $\mathscr C_s$ is linear and continuous, there is $K<\infty$ such that
\begin{equation*}
    \big\|\mathscr C_s(\widetilde X^k-\widetilde X^\star)\big\|_2<\delta
    \qquad\text{for every }k\geq K.
\end{equation*}
The uniqueness just proved forces every exact update to satisfy
$S^{k+1}=S^\star$.  In fact, locking is permanent after the first such update.
Indeed, define
\begin{equation*}
    y^\star=(\mathcal A\mathcal A^*)^{-1}\mathcal A(S^\star+C).
\end{equation*}
The feasibility 
of $S^\star$ implies
$\mathcal A^*(y^\star)=\mathcal P(S^\star+C)=S^\star+C$, so the residual is
zero and the multiplier update gives $\widetilde X^{k+1}=\widetilde X^k$.
The same certificate $X$ consequently remains valid for every subsequent
value of $\sigma>0$, and induction proves \eqref{eq:finite-locking-conclusion}.
Finally,
\begin{equation*}
    \mathcal A(X+\Diag(z))
    =\mathcal A(\widetilde X+D)=b,
    \qquad
    X\succeq0,
    \qquad
    XS^\star=0.
\end{equation*}
Hence the locked tuple satisfies the full
KKT system of \eqref{dsdp}.
\end{proof}


The neighborhood condition in Theorem~\ref{thm:finite-locking} can be replaced
by a simpler, although more conservative, condition.  Indeed, for every
$G\in\mathbb S_n$,
\begin{align}
    \|\mathscr C_s(G)\|_2\leq\|G\|_2
    +\|\Diag((Gs)\circ s)\|_2=\|G\|_2+\|Gs\|_\infty
    \leq(1+\sqrt n)\|G\|_2.                       \label{eq:certificate-lipschitz}
\end{align}
Consequently,
\begin{equation}\label{eq:simple-locking-neighborhood}
    \|\widetilde X-\widetilde X^\star\|_2
    <\frac{\delta}{1+\sqrt n}
\end{equation}
implies \eqref{eq:locking-neighborhood}.  If, for example,
$\|\widetilde X^k-\widetilde X^\star\|_2\leq c\rho^k$ for some $c>0$ and
$\rho\in(0,1)$, then locking occurs no later than any integer $K$ for which
\begin{equation}\label{eq:locking-iteration-estimate}
    (1+\sqrt n)c\rho^K<\delta.
\end{equation}

\begin{remark}\label{rem:inexact-near-locking}
Literal finite locking does not follow solely from the inexact subproblem stopping criterion. Nevertheless, Theorem~\ref{thm:finite-locking}
gives a quantitative near-locking statement.  If
$\widehat S\in\mathcal M$ satisfies
\begin{equation*}
    \Phi_{\widetilde X,\sigma}(\widehat S)
    \leq
    \min_{S\in\mathcal M}\Phi_{\widetilde X,\sigma}(S)+\eta,
\end{equation*}
where \eqref{eq:locking-neighborhood} holds, then
\begin{equation}\label{eq:inexact-near-locking}
    \|\widehat S-S^\star\|
    \leq\frac{n\eta}{\gamma(\widetilde X)},
    \qquad
    \|\widehat S-S^\star\|
    \leq\sqrt{\frac{2n\eta}{\gamma(\widetilde X)}}.
\end{equation}
Thus, solving the Riemannian subproblem with increasing accuracies drives the computed matrix
sharply toward the rank-one solution. An exact rank-compression or rounding
step can convert this near-locking into literal finite locking, which provides
a mechanism for the residue diving phenomenon observed in practice.
\end{remark}

\section{Numerical experiments}\label{experiments}
We implemented \texttt{ManiDSDP} in MATLAB within the \texttt{ManiSDP} codebase.\footnote{The codebase is available at \url{https://github.com/wangjie212/ManiSDP-matlab}.} The factorized ALM subproblems are solved with \texttt{Manopt 7.0} \cite{manopt}. All reported times are wall-clock times in seconds. A dash indicates that the solver ran out of memory, $*$ indicates that the running time exceeded $10{,}000$ seconds, and $**$ indicates a numerical failure.

{\bf Hardware.} All experiments were performed on a desktop computer with an Intel Core i9-10900 processor at $2.80$~GHz and $64$~GB of RAM.

{\bf Baseline Solvers.} We compare \texttt{ManiDSDP} with three complementary baselines: \texttt{MOSEK 11.0} \cite{mosek}, a widely used commercial conic solver; \texttt{SDPNAL+} \cite{sun2020sdpnal}, which targets large-scale SDPs; and \texttt{ManiSDP} \cite{wang2025solving}, a recent primal low-rank Riemannian ALM solver.

{\bf The Stopping Criterion.} For an approximate primal-dual tuple $(S,y,X,z)$, define
\begin{align}\label{eq:res}
\eta_p &= \frac{\left\lVert\cA^*(y)-S-C\right\rVert}{1 +\left\lVert C \right\rVert},\\
\eta_d &= \frac{\max\{0,-\lambda_{\min}(X)\}}{1+|\lambda_{\max}(X)|},\\
\eta_g &= \frac{|\langle C, X+\cB^*(z)\rangle+\sum_{i=1}^nz_i- b^{\intercal}y|}{1 + |\langle C, X+\cB^*(z)\rangle+\sum_{i=1}^nz_i| + |b^{\intercal}y|}.
\end{align}
We report $\eta_{\max}\coloneqq\max\{\eta_p,\eta_d,\eta_g\}$ and terminate when $\eta_{\max}\leq10^{-8}$.

\subsection{Dense BQPs}
We first consider the dense BQP \eqref{sec3:bqp}, a model that includes Max-Cut, $\mathbb{Z}_2$ synchronization, and the classical Ising problem.
For each $q\in\{10,20,30,40,50,60,70,80\}$, we generate three random instances by drawing the independent entries of $Q\in\S_q$ and the entries of $\mathbf c\in\R^q$ from $\mathcal N(0,1)$. For each instance, we solve the second-order moment-SOS relaxation described in Section~\eqref{sec3-1}. We set $p_0=\lceil\log(m)\rceil$. The results are presented in Table~\ref{table:bpq1}.

\begin{table}[htbp]
\caption{Results for dense BQPs.}\label{table:bpq1}
\renewcommand\arraystretch{1.2}
\centering
\begin{tabular}{c|c|c|c|c|c|c|c|c|c}
\multirow{2}{*}{$q$}&\multirow{2}{*}{trial}&\multicolumn{2}{c|}{{\tt MOSEK}}&\multicolumn{2}{c|}{{\tt SDPNAL+}}&\multicolumn{2}{c|}{{\tt ManiSDP}}&\multicolumn{2}{c}{{\tt ManiDSDP}}\\
		\cline{3-10}
&&$\eta_{\max}$&time&$\eta_{\max}$&time&$\eta_{\max}$&time&$\eta_{\max}$&time\\
		\hline
		\multirow{3}{*}{10}&\#1&2.8e-09&{\bf 0.02}&1.5e-09&0.10&3.9e-15&0.07&8.4e-15&0.07\\
        &\#2&7.4e-10&{\bf 0.02}&4.4e-10&0.08&4.1e-15&0.06&8.8e-15&0.06\\
        &\#3&3.0e-11&{\bf 0.02}&6.0e-09&0.08&3.7e-15&0.07&8.9e-15&0.07\\
		\hline
		\multirow{3}{*}{20}&\#1&3.1e-13&3.49&3.8e-10&1.21&1.6e-14&0.25&4.9e-14&{\bf 0.20}\\
        &\#2&1.0e-09&3.64&3.1e-09&1.09&1.3e-14&0.33&4.8e-14&{\bf 0.24}\\
        &\#3&1.6e-10&3.44&2.0e-09&1.46&1.4e-14&0.42&4.8e-14&{\bf 0.28}\\
		\hline
		\multirow{3}{*}{30}&\#1&1.3e-12&439&9.4e-10&6.74&2.8e-14&3.14&1.1e-13&{\bf 1.64}\\
        &\#2&5.2e-12&426&7.2e-09&14.2&2.9e-14&1.67&1.1e-13&{\bf 1.15}\\
        &\#3&1.8e-10&491&3.4e-10&20.9&3.0e-14&2.57&1.1e-13&{\bf 1.52}\\
		\hline
		\multirow{3}{*}{40}&\#1&-&-&2.4e-09&79.2&4.9e-14&7.93&2.3e-13&{\bf 3.36}\\
        &\#2&-&-&4.1e-10&92.6&4.6e-14&4.83&2.4e-13&{\bf 2.01}\\
        &\#3&-&-&5.0e-10&65.6&4.8e-14&7.71&2.3e-13&{\bf 2.71}\\
		\hline
		\multirow{3}{*}{50}&\#1&-&-&6.7e-08&558&7.0e-14&25.9&4.3e-13&{\bf 11.1}\\
        &\#2&-&-&6.1e-11&480&2.9e-13&34.2&4.3e-13&{\bf 14.3}\\
        &\#3&-&-&3.8e-09&607&6.8e-14&58.1&4.3e-13&{\bf 19.6}\\
		\hline
		\multirow{3}{*}{60}&\#1&-&-&2.8e-07&1832&9.6e-14&80.6&7.5e-13&{\bf 58.6}\\
        &\#2&-&-&9.1e-05&7464&9.7e-14&597&7.5e-13&{\bf 234}\\
        &\#3&-&-&1.6e-07&1898&9.6e-14&146&7.5e-13&{\bf 53.7}\\
        \hline
        \multirow{3}{*}{70}&\#1&-&-&$*$&$*$&1.4e-13&986&7.7e-13&{\bf 366}\\
        &\#2&-&-&9.8e-07&6934&1.4e-13&272&7.7e-13&{\bf 142}\\
        &\#3&-&-&$*$&$*$&1.2e-13&421&7.7e-13&{\bf 172}\\
        \hline
        \multirow{3}{*}{80}&\#1&-&-&$*$&$*$&1.6e-13&1037&1.3e-12&{\bf 492}\\
        &\#2&-&-&$*$&$*$&1.7e-13&1844&2.5e-11&{\bf 656}\\
        &\#3&-&-&$*$&$*$&1.6e-13&2145&3.1e-12&{\bf 1181}\\
	\end{tabular}
\end{table}

Table~\ref{table:bpq1} reveals a clear change in the relative performance of the solvers as the problem size grows. For $q=10$, \texttt{MOSEK} is the fastest solver, but its running time increases sharply and it runs out of memory for $q\geq 40$. \texttt{SDPNAL+} becomes increasingly expensive and cannot complete most instances with $q\geq 70$. In contrast, both low-rank SDP solvers complete all tested instances. Among them, \texttt{ManiDSDP} is the fastest solver in every trial with $q\geq 20$, while maintaining KKT residues extremely small. Finally, we mention that in these experiments,, the final factorization size is $2$ and the output matrices $S,X$ of \texttt{ManiDSDP} satisfy
\begin{equation}
    \rank(S)=1 \qquad \rank(X)=n-1.
\end{equation}
Thus, the observed solutions are consistent with both the low-rank premise and strict complementarity.




\subsection{Sparse BQPs}
We next consider sparse BQPs with overlapping variable blocks. Let $t,q\geq1$ be integers. For $k=1,\ldots,t$, define
\[
    \x_k\coloneqq
    \bigl(x_{(q-2)(k-1)+1},\ldots,x_{(q-2)k+2}\bigr)^{\intercal}\in\R^q,
\]
and let $\x=(x_1,\ldots,x_{(q-2)t+2})^{\intercal}$. The resulting problem is
\begin{equation}\label{sec6:sbqp}
\begin{cases}
\inf\limits_{\x\in\R^{(q-2)t+2}} & \sum_{k=1}^t\left(\x_k^\intercal Q_k\x_k + \mathbf{c}_k^\intercal\x_k\right)\\
\,\,\,\,\quad\rm{s.t.}&x_i^2=1,\quad i=1,2,\ldots,(q-2)t+2,\\
\end{cases}\tag{BQP-sparse}
\end{equation}
where $Q_k\in\S_q$ and $\mathbf{c}_k\in\R^{q}$, $k=1,2,\ldots,t$.
The second-order sparse SOS relaxation of \eqref{sec6:sbqp} is the following multi-block SDP:
\begin{equation}\label{sec6:ssos}
\begin{cases}
\sup\limits_{X_k,\lambda} &\lambda\\
\,\,\rm{s.t.}&\sum_{k=1}^t(\x_k^\intercal Q_k\x_k + \mathbf{c}_k^\intercal\x_k)-\lambda\equiv \sum_{k=1}^tv(\x_k)^{\intercal}X_kv(\x_k)\quad\mathrm{mod}\,\,I,\\
&X_k\succeq0,\quad k=1,2,\ldots,t,\\
\end{cases}
\end{equation}
where $I$ is the ideal generated by $\{x_i^2-1\}_{i=1}^{(q-2)t+2}$. The dual formulation of \eqref{sec6:ssos} fits into a multi-block version of \eqref{dsdp} with $C=0$:
\begin{equation}\label{sdsdp}
\begin{cases}
\inf\limits_{y\in\R^m,S\in\S_n}&b^{\intercal}y\\
\,\,\,\,\quad\rm{s.t.}&S=\cA^*(y)-C\in\S^+_{|v(\x_1)|}\times\cdots\times\S^+_{|v(\x_t)|},\\
&\diag(S)=\mathbf{1}.
\end{cases}\tag{DSDP-sparse}
\end{equation}

First, we fix the parameter $q=20$ and vary the number of blocks. For each $t\in\{20,40,60,80,100,120\}$, we generate three random instances of \eqref{sec6:sbqp} by taking $Q_k\in\S_q$ with $[Q_k]_{ij}\sim\cN(0,1)$ and $\mathbf{c}_k\in\R^{q}$ with $c_{k,i}\sim\cN(0,1)$, $k=1,2,\ldots,t$. For each instance, we solve the sparse SDP relaxation \eqref{sdsdp}. The results are presented in Table~\ref{table:bpq2}. 

\begin{table}[htbp]
	\caption{Results for sparse BQPs ($q=20$).}\label{table:bpq2}
	\centering
	\begin{tabular}{c|c|c|c|c|c|c|c|c|c}
		\multirow{2}{*}{$t$}&\multirow{2}{*}{trial}&\multicolumn{2}{c|}{{\tt MOSEK}}&\multicolumn{2}{c|}{{\tt SDPNAL+}}&\multicolumn{2}{c|}{{\tt ManiSDP}}&\multicolumn{2}{c}{{\tt ManiDSDP}}\\
		\cline{3-10}
        &&$\eta_{\max}$&time&$\eta_{\max}$&time&$\eta_{\max}$&time&$\eta_{\max}$&time\\
		\hline
		\multirow{3}{*}{20}&\#1&1.4e-13&151&5.5e-10&226&1.9e-17&26.0&1.3e-12&{\bf 12.5}\\
             &\#2&6.3e-14&146&1.6e-08&143&4.7e-14&20.4&1.3e-12&{\bf 6.97}\\
             &\#3&1.2e-13&160&6.7e-08&275&1.4e-17&34.3&1.3e-12&{\bf 11.4}\\
             \hline
            \multirow{3}{*}{40}&\#1&8.1e-14&283&3.5e-09&622&4.1e-14&77.5&7.0e-12&{\bf 51.0}\\
            &\#2&8.2e-14&340&4.2e-08&1072&1.9e-17&130&7.0e-12&{\bf 45.1}\\
            &\#3&1.2e-13&289&7.1e-08&1206&1.9e-17&237&7.0e-12&{\bf 49.2}\\
            \hline
            \multirow{3}{*}{60}&\#1&1.4e-13&400&4.9e-10&2206&5.7e-13&79.7&1.5e-11&{\bf 67.9}\\
            &\#2&2.2e-13&448&4.5e-07&3322&3.9e-14&313&1.5e-11&{\bf 80.3}\\
            &\#3&1.5e-13&507&6.7e-10&1449&9.8e-14&167&1.5e-11&{\bf 154}\\
            \hline
            \multirow{3}{*}{80}&\#1&1.7e-12&651&7.2e-10&2933&9.3e-17&824&1.6e-11&{\bf 304}\\
            &\#2&2.4e-13&718&2.7e-09&4542&9.3e-17&342&1.6e-11&{\bf 119}\\
            &\#3&1.0e-13&719&2.0e-09&3195&1.8e-16&153&1.6e-11&{\bf 144}\\
            \hline
            \multirow{3}{*}{100}&\#1&2.3e-12&633&3.4e-09&3581&7.4e-17&705&2.5e-11&{\bf 287}\\
            &\#2&5.0e-13&866&1.9e-08&6811&1.5e-16&1595&2.5e-11&{\bf 80.9}\\
            &\#3&5.9e-13&956&1.6e-09&3249&7.0e-14&231&2.5e-11&{\bf 131}\\
            \hline
            \multirow{3}{*}{120}&\#1&4.4e-13&1269&1.8e-10&6931&1.2e-16&1175&3.4e-11&{\bf 573}\\
            &\#2&2.5e-13&1039&8.6e-11&4942&1.1e-13&962&3.4e-11&{\bf 135}\\
            &\#3&8.4e-13&910&5.0e-09&4863&6.1e-17&278&3.4e-11&{\bf 228}\\
	\end{tabular}
\end{table}

Table~\ref{table:bpq2} studies scalability with respect to the number of PSD blocks while keeping the block size fixed. \texttt{ManiDSDP} solves every instance with high accuracy and is the fastest solver in every trial. \texttt{MOSEK} also remains accurate, but is consistently more expensive. The running time of \texttt{SDPNAL+} grows rapidly and some of its reported residues are above the prescribed tolerance. \texttt{ManiSDP} generally attains very small residues, but its running time varies considerably across instances. These results indicate that \texttt{ManiDSDP} exploits the multi-block sparsity effectively and remains robust as the number of blocks increases.

Next, we fix $t=10$ and vary the block size. For each $q\in\{10,20,30,40\}$, we generate three random instances of \eqref{sec6:sbqp} by taking $Q_k\in\S_q$ with $[Q_k]_{ij}\sim\cN(0,1)$ and $\mathbf{c}_k\in\R^{q}$ with $c_{k,i}\sim\cN(0,1)$, $k=1,2,\ldots,t$. For each instance, we solve the sparse SDP relaxation \eqref{sdsdp}. The results are presented in Table~\ref{table:bpq3}. 

\begin{table}[htbp]
	\caption{Results for sparse BQPs ($t=10$).}\label{table:bpq3}
	\centering
	\begin{tabular}{c|c|c|c|c|c|c|c|c|c}
		\multirow{2}{*}{$q$}&\multirow{2}{*}{trial}&\multicolumn{2}{c|}{{\tt MOSEK}}&\multicolumn{2}{c|}{{\tt SDPNAL+}}&\multicolumn{2}{c|}{{\tt ManiSDP}}&\multicolumn{2}{c}{{\tt ManiDSDP}}\\
		\cline{3-10}
&&$\eta_{\max}$&time&$\eta_{\max}$&time&$\eta_{\max}$&time&$\eta_{\max}$&time\\
		\hline
		\multirow{3}{*}{10}
             &\#1&1.2e-12&{\bf 0.13}&1.0e-09&0.75&4.5e-14&0.31&4.3e-14&0.14\\
             &\#2&2.3e-11&{\bf 0.13}&7.0e-09&0.87&6.7e-17&0.23&3.3e-13&0.19\\
             &\#3&5.2e-13&{\bf 0.14}&4.7e-09&1.04&5.5e-17&0.34&2.1e-11&0.21\\
             \hline
		\multirow{3}{*}{20}&\#1&3.2e-13&58.6&5.7e-08&90.6&2.1e-17&13.4&1.2e-12&{\bf 4.37}\\
             &\#2&1.3e-13&64.1&4.4e-09&39.9&8.8e-17&11.7&1.2e-12&{\bf 3.34}\\
             &\#3&9.4e-14&60.0&1.7e-09&67.2&2.2e-17&19.6&1.6e-12&{\bf 4.93}\\
             \hline
            \multirow{3}{*}{30}&\#1&-&-&9.1e-09&1092&1.1e-17&145&1.7e-12&{\bf 58.0}\\
            &\#2&-&-&5.3e-09&619&9.6e-17&250&1.7e-12&{\bf 59.4}\\
            &\#3&-&-&9.3e-09&1087&8.4e-18&208&1.7e-12&{\bf 71.2}\\
            \hline
            \multirow{3}{*}{40}&\#1&-&-&$*$&$*$&6.4e-17&3596&7.9e-12&{\bf 756}\\
            &\#2&-&-&4.0e-09&3287&1.8e-16&1639&7.9e-12&{\bf 566}\\
            &\#3&-&-&2.9e-07&6480&$**$&$**$&7.9e-12&{\bf 1114}\\
	\end{tabular}
\end{table}

Table~\ref{table:bpq3} complements the preceding experiment by fixing $t=10$ and increasing the block size. For $q=10$, \texttt{MOSEK} is marginally the fastest solver. For every trial with $q\geq 20$, however, \texttt{ManiDSDP} has the shortest running time. Moreover, it is the only solver that solves all three instances with $q=40$ to the required tolerance: \texttt{MOSEK} runs out of memory for $q\geq 30$, and both \texttt{SDPNAL+} and \texttt{ManiSDP} fail to complete one $q=40$ instance. 

\subsection{UCQPs}
We next test the real-coefficient UCQP \eqref{eq:ucqp-size-comparison}, a model that includes formulations arising in phase recovery, unimodular code design, and multiple-input multiple-output detection. For each $q\in\{10,20,30,40\}$, we generate three random instances by drawing the independent entries of $Q\in\S_q$ and the entries of $\mathbf c\in\R^q$ from $\mathcal N(0,1)$. For each instance, we solve the second-order moment-SOS relaxation described in Section~\ref{sec3-1}. We set $p_0=\lceil\log(m)\rceil$. Table~\ref{table:uqp} reports the results.

\begin{table}[htbp]
	\caption{Results for UCQPs.}\label{table:uqp}
	\centering
	\begin{tabular}{c|c|c|c|c|c|c|c|c|c}
		\multirow{2}{*}{$q$}&\multirow{2}{*}{trial}&\multicolumn{2}{c|}{{\tt MOSEK}}&\multicolumn{2}{c|}{{\tt SDPNAL+}}&\multicolumn{2}{c|}{{\tt ManiSDP}}&\multicolumn{2}{c}{{\tt ManiDSDP}}\\
		\cline{3-10}
&&$\eta_{\max}$&time&$\eta_{\max}$&time&$\eta_{\max}$&time&$\eta_{\max}$&time\\
		\hline
        \multirow{3}{*}{10}&\#1&7.2e-09&1.76&7.8e-10&1.39&2.4e-09&0.35&3.9e-09&{\bf 0.19}\\
            &\#2&1.0e-10&2.10&7.5e-09&2.88&8.8e-09&0.33&8.5e-09&{\bf 0.19}\\
            &\#3&2.6e-12&2.17&5.1e-09&0.85&4.7e-10&0.26&8.2e-09&{\bf 0.16}\\
            \hline
            \multirow{3}{*}{20}&\#1&4.6e-11&3808&5.6e-09&216&5.6e-10&11.5&2.9e-09&{\bf 4.72}\\
            &\#2&9.5e-12&3523&7.8e-07&147&6.9e-09&5.91&6.7e-09&{\bf 5.04}\\
            &\#3&1.5e-10&3552&9.9e-07&137&5.1e-09&8.20&2.6e-09&{\bf 3.77}\\
            \hline
            \multirow{3}{*}{30}&\#1&-&-&9.2e-11&1931&7.3e-10&210&7.5e-09&{\bf 53.6}\\
            &\#2&-&-&6.8e-09&1355&3.5e-09&192&2.2e-09&{\bf 34.0}\\
            &\#3&-&-&1.3e-05&6973&3.4e-09&120&9.2e-09&{\bf 45.2}\\
            \hline
            \multirow{3}{*}{40}&\#1&-&-&$*$&$*$&8.8e-09&{\bf 1030}&3.4e-09&1438\\
            &\#2&-&-&$*$&$*$&2.0e-09&1140&6.7e-09&{\bf 309}\\
            &\#3&-&-&$*$&$*$&$**$&$**$&9.5e-09&{\bf 481}\\
	\end{tabular}
\end{table}

Table~\ref{table:uqp} shows that the advantage of the proposed method also extends to SDP relaxations of complex polynomial optimization problems. \texttt{ManiDSDP} meets the prescribed stopping tolerance on all instances and is the fastest solver in every trial for $q\leq 30$. At $q=40$, \texttt{ManiSDP} is faster on the first trial, but \texttt{ManiDSDP} is faster on the other two trials and is the only solver that solves all three instances to the requested accuracy. In comparison, \texttt{MOSEK} runs out of memory from $q=30$ onward, \texttt{SDPNAL+} exceeds the time limit at $q=40$, and \texttt{ManiSDP} encounters some numerical issue on one of the largest instances. Hence, the principal benefit of \texttt{ManiDSDP} on UCQPs is not only speed but also robustness at the largest tested dimension. Finally, we mention that for these tests, the final factorization size is $2$ and the output matrices $S,X$ of \texttt{ManiDSDP} satisfy
\begin{equation}
    \rank(S)=2 \qquad \rank(X)=n-2.
\end{equation}
Thus, the observed solutions are therefore consistent with the low-rank premise and strict complementarity.

\subsection{{\tt ManiDSDP} versus {\tt ManiSDP}}
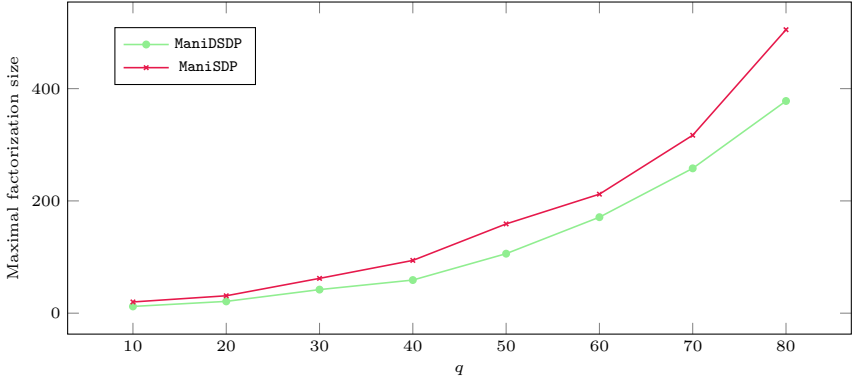
\begin{figure}[htbp]	
\centering
\begin{tikzpicture}
\footnotesize
\scalefont{0.8} 
\begin{axis}[
sharp plot, 
xmode=normal,
width=12cm, height=6cm,  
xlabel = $q$,
ylabel = Maximal factorization size,
xlabel near ticks,
ylabel near ticks,
legend style={at={(0.15,0.75)},anchor=south},
]

\addplot[semithick,mark=*,mark options={scale=0.6}, color=lightgreen] coordinates { 
    (10,12)
    (20,21)
    (30,42)
    (40,59)
    (50,106)
    (60,171)
    (70,258)
    (80,378)
  };
\addlegendentry{{\tt ManiDSDP}}

\addplot[semithick,mark=x,mark options={scale=0.6}, color=color3] coordinates { 
    (10,20)
    (20,31)
    (30,62)
    (40,94)
    (50,159)
    (60,212)
    (70,317)
    (80,505)
};
\addlegendentry{{\tt ManiSDP}}
\end{axis}
\end{tikzpicture}
\caption{Comparison of maximal factorization sizes reached through outer ALM iterations of {\tt ManiDSDP} and {\tt ManiSDP}.}
\label{fig:0}
\end{figure}

\begin{figure}[htbp]	
\centering
\begin{tikzpicture}
\footnotesize
\scalefont{0.8} 
\begin{axis}[
sharp plot, 
xmode=normal,
width=12cm, height=6cm,  
xlabel = $q$,
ylabel = Number of outer ALM iterations,
xlabel near ticks,
ylabel near ticks,
legend style={at={(0.15,0.75)},anchor=south},
]

\addplot[semithick,mark=*,mark options={scale=0.6}, color=lightgreen] coordinates { 
    (10,13)
    (20,13)
    (30,17)
    (40,12)
    (50,20)
    (60,50)
    (70,56)
    (80,75)
  };
\addlegendentry{{\tt ManiDSDP}}

\addplot[semithick,mark=x,mark options={scale=0.6}, color=color3] coordinates { 
    (10,13)
    (20,16)
    (30,27)
    (40,22)
    (50,37)
    (60,116)
    (70,121)
    (80,122)
};
\addlegendentry{{\tt ManiSDP}}
\end{axis}
\end{tikzpicture}
\caption{Comparison of numbers of outer ALM iterations taken by {\tt ManiDSDP} and {\tt ManiSDP}.}
\label{fig:4}
\end{figure}
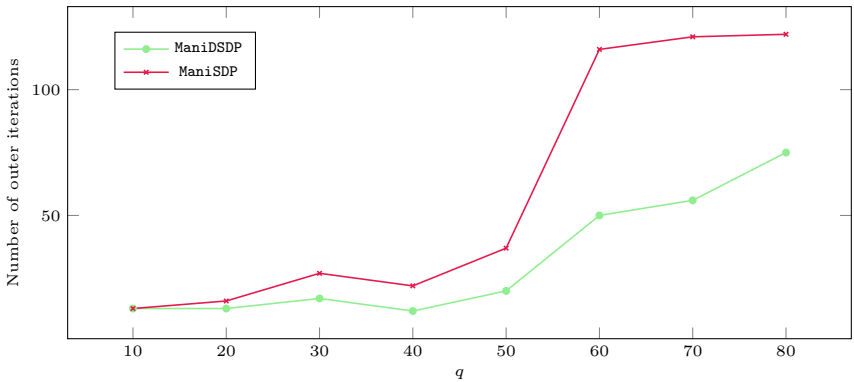

In this subsection, let us compare {\tt ManiDSDP} with {\tt ManiSDP} in more detail. 
In Fig.~\ref{fig:0}, we display the maximal factorization size reached through outer ALM iterations of {\tt ManiDSDP} and {\tt ManiSDP} (each point is averaged over three random instances of \eqref{sec3:bqp}). 
In Fig.~\ref{fig:4}, we display the number of outer ALM iterations taken by {\tt ManiDSDP} and {\tt ManiSDP} (each point is averaged over three random instances of \eqref{sec3:bqp}).

Figures~\ref{fig:0} and~\ref{fig:4} identify two complementary sources of the runtime improvement of \texttt{ManiDSDP} over \texttt{ManiSDP}. The maximal factorization size used by \texttt{ManiDSDP} is uniformly smaller, amounting to roughly $60\%$--$75\%$ of that used by \texttt{ManiSDP} over the tested range. The two methods take the same number of outer ALM iterations at $q=10$, but \texttt{ManiDSDP} requires fewer iterations for every $q\geq 20$, with the gap becoming particularly pronounced for $q=60,70,80$. Since the cost of the factorized subproblems grows rapidly with the factorization size, the simultaneous reduction in both factorization sizes and numbers of outer ALM iterations explains a substantial part of the timing advantage observed in Table~\ref{table:bpq1}.

\subsection{\eqref{dsdp1} versus \eqref{dsdp}}\label{sec:reformulation-comparison-experiments}
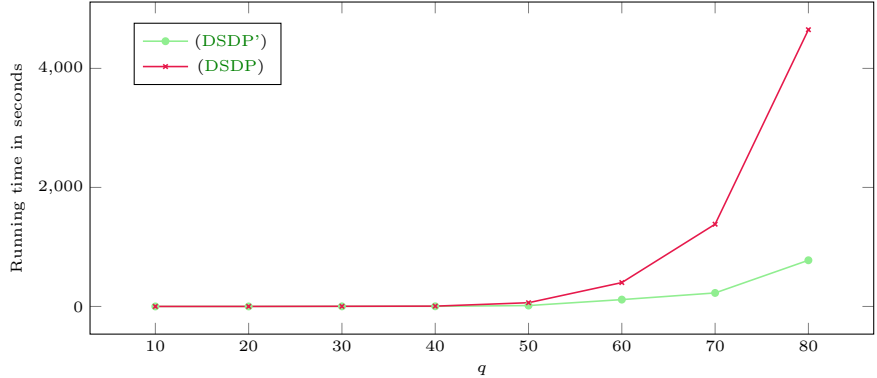
\begin{figure}[htbp]	
\centering
\begin{tikzpicture}
\footnotesize
\scalefont{0.8} 
\begin{axis}[
sharp plot, 
xmode=normal,
width=12cm, height=6cm,  
xlabel = $q$,
ylabel = Running time in seconds,
xlabel near ticks,
ylabel near ticks,
legend style={at={(0.15,0.75)},anchor=south},
]

\addplot[semithick,mark=*,mark options={scale=0.6}, color=lightgreen] coordinates { 
    (10,0.07)
    (20,0.24)
    (30,1.44)
    (40,2.69)
    (50,15.0)
    (60,115)
    (70,227)
    (80,776)
  };
\addlegendentry{\eqref{dsdp1}}

\addplot[semithick,mark=x,mark options={scale=0.6}, color=color3] coordinates { 
    (10,0.11)
    (20,0.40)
    (30,2.70)
    (40,5.32)
    (50,62.6)
    (60,401)
    (70,1382)
    (80,4650)
};
\addlegendentry{\eqref{dsdp}}
\end{axis}
\end{tikzpicture}
\caption{Comparison of running time in solving \eqref{dsdp1} and \eqref{dsdp}.}
\label{fig:1}
\end{figure}

In this subsection, we compare the running time in solving \eqref{dsdp1} and \eqref{dsdp} with the dual Riemannian ALM/ADMM algorithm, respectively. We show the results in Fig.~\ref{fig:1} (each point is averaged over three random instances of \eqref{sec3:bqp}). 

Figure~\ref{fig:1} shows that the variable-projected formulation \eqref{dsdp1} is faster than the direct formulation \eqref{dsdp} for every tested value of $q$. The speedup is about $1.6$ for $q=10$ and grows to approximately $6$ for $q=70$ and $q=80$. This increasing separation indicates that eliminating the variable $y$ and imposing affine consistency through the orthogonal-projector residual becomes progressively more beneficial as the problem dimension increases. The experiment therefore provides direct empirical justification for using \eqref{dsdp1} as the computational formulation of the dual Riemannian ALM.

\subsection{The residue diving phenomenon}\label{rdp}
On random BQP instances, \texttt{ManiDSDP} exhibits a \emph{residue diving phenomenon}: at a terminal outer iteration, the maximum KKT residual drops abruptly to far below $10^{-8}$. Figure~\ref{fig:3} illustrates this behavior for one representative random instance of \eqref{sec3:bqp} at each $q\in\{10,20,30,40,50,60\}$. Similar behavior is also observed with \texttt{ManiSDP}.

\begin{figure}[htbp]	
\centering
\begin{tikzpicture}
\footnotesize
\scalefont{0.8} 
\begin{axis}[
sharp plot, 
xmode=normal,
width=12cm, height=6cm,  
xlabel = Iteration,
ylabel = $\log_{10}\eta_{\max}$,
xlabel near ticks,
ylabel near ticks,
legend style={at={(0.15,0.05)},anchor=south},
]

\addplot[semithick,mark=*,mark options={scale=0.6}, color=lightgreen] coordinates { 
   (1,0.5965)
   (2,0.5039)
    (3,0.4318)
    (4,0.3372)
    (5,0.1652)
   (6,-0.1739)
   (7,-0.9046)
   (8,-1.9112)
   (9,-2.9540)
   (10,-3.5970)
   (11,-4.2637)
   (12,-4.9166)
  (13,-14.0753)
  };
\addlegendentry{$q=10$}

\addplot[semithick,mark=x,mark options={scale=0.6}, color=bordeaux] coordinates { 
    (1,1.1509)
    (2,0.9821)
    (3,0.8113)
    (4,0.6613)
    (5,0.4579)
    (6,0.2072)
   (7,-0.1369)
   (8,-0.6413)
   (9,-1.7460)
   (10,-2.7107)
   (11,-3.5207)
   (12,-4.2325)
  (13,-13.3129)
};
\addlegendentry{$q=20$}

 \addplot[semithick,mark=+,mark options={scale=0.6}, color=color1] coordinates { 
    (1,1.4752)
    (2,1.3162)
    (3,1.0354)
    (4,0.8495)
    (5,0.6570)
    (6,0.3128)
   (7,-0.0985)
   (8,-0.4180)
   (9,-0.7281)
   (10,-1.0094)
   (11,-1.2710)
   (12,-1.5707)
   (13,-1.7545)
   (14,-2.0256)
   (15,-1.8494)
   (16,-2.5981)
   (17,-3.7050)
   (18,-4.4020)
  (19,-12.9443)
};
\addlegendentry{$q=30$}

 \addplot[semithick,mark=*,mark options={scale=0.6}, color=color2] coordinates { 
    (1,1.8019)
    (2,1.5462)
    (3,1.4162)
    (4,1.2220)
    (5,1.0201)
    (6,0.5679)
    (7,0.2985)
   (8,-0.0729)
   (9,-0.6385)
   (10,-1.2715)
   (11,-2.1721)
   (12,-3.0854)
   (13,-4.0871)
  (14,-12.6290)

};
    \addlegendentry{$q=40$}

\addplot[semithick,mark=x,mark options={scale=0.6}, color=color3] coordinates {
    (1,1.9666)
    (2,1.8333)
    (3,1.7055)
    (4,1.2557)
    (5,0.8682)
    (6,0.6333)
    (7,0.3832)
    (8,0.1435)
   (9,-0.0912)
   (10,-0.2978)
   (11,-0.3938)
   (12,-0.9625)
   (13,-1.7298)
   (14,-2.6869)
   (15,-3.5540)
  (16,-12.3644)
};
    \addlegendentry{$q=50$}

    \addplot[semithick,mark=+,mark options={scale=0.6}, color=darkblue] coordinates { 
     (1,2.0757)
    (2,1.9840)
    (3,1.6189)
    (4,1.3956)
    (5,1.3079)
    (6,1.1906)
    (7,1.1029)
    (8,1.0595)
    (9,1.0121)
    (10,0.9321)
    (11,0.7994)
    (12,0.5990)
    (13,0.4966)
    (14,0.4343)
    (15,0.4273)
   (16,-0.0688)
    (17,0.0317)
    (18,0.1949)
   (19,-0.0386)
   (20,-0.0898)
   (21,-0.2058)
   (22,-0.0841)
   (23,-0.1456)
   (24,-0.4082)
   (25,-0.4938)
   (26,-0.7630)
   (27,-1.0306)
   (28,-1.1005)
   (29,-1.0104)
   (30,-1.1443)
   (31,-1.0776)
   (32,-1.1546)
   (33,-1.2411)
   (34,-1.5219)
   (35,-2.0322)
   (36,-3.4020)
  (37,-12.1260)
};
    \addlegendentry{$q=60$}
    
\end{axis}
\end{tikzpicture}
\caption{Maximum KKT residual at each outer iteration for random BQPs solved by \texttt{ManiDSDP}.}
\label{fig:3}
\end{figure}  

The curves in Figure~\ref{fig:3} exhibit a common two-phase pattern. During the first phase, the maximal KKT residue decreases gradually, and for the larger instances the decrease can be nonmonotone or contain a prolonged plateau. Once the iterations enter the terminal regime, however, the residue drops in a single outer iteration from approximately $10^{-3}$--$10^{-5}$ to $10^{-12}$--$10^{-14}$. This observation is consistent with the local one-step locking result proved in Section~\ref{finite-locking}.

\section{Conclusions and future work}\label{conclusions}
We have developed \texttt{ManiDSDP}, a dual Riemannian augmented Lagrangian method for low-rank SDPs with unit diagonal constraints. The method combines three structural reductions. First, under the invertibility of $\cA\cA^*$, the variable $y$ is eliminated from SDPs by a variable projection, and affine consistency is represented by an orthogonal-projector residual. Second, the factorization $S=YY^{\intercal}$ converts the unit-diagonal condition into an oblique-manifold constraint. Third, dynamic rank adjustment and negative-curvature information are used to avoid terminating at nonoptimal stationary points of the factorized subproblems.

The analysis establishes several convergence and identification properties. Under controlled subproblem tolerances, every cluster point is a KKT point. Under the additional assumptions of full-sequence convergence, strict complementarity, and finite-tail exact complementarity, the rank of the limiting solution is identified after finitely many iterations. Under a strictly complementary rank-one solution, the local one-step locking theorem further shows that, once the multiplier enters an explicit neighborhood, the exact subproblem recovers the rank-one solution uniquely and the subsequent iterations remain there. This result gives a local theoretical explanation for the terminal residue dive observed in the experiments.

The numerical results on dense and sparse BQPs and on UCQPs indicate that \texttt{ManiDSDP} can attain high accuracy while remaining efficient on instances for which conventional SDP solvers become expensive or fail because of time or memory limitations. The comparison with \texttt{ManiSDP} shows that the dual formulation typically uses smaller factorization sizes and fewer outer iterations on the tested problems. Moreover, the direct comparison between \eqref{dsdp1} and \eqref{dsdp} confirms that the variable projection yields an increasingly significant computational advantage as the problem dimension grows. 

There are several directions for future research.
\begin{itemize}

    \item \emph{Convergence rates and iteration complexity.}
    The local one-step locking theorem explains the terminal phase under a strictly complementary rank-one solution, but does not describe the preceding phase of the algorithm. Establishing global or local rates under error-bound, metric-subregularity, or Kurdyka--{\L}ojasiewicz-type conditions would clarify when the observed rapid convergence can be expected.
    

    \item \emph{Weaker structural assumptions.}
    The current variable projection uses the invertibility of $\cA\cA^*$, and the manifold formulation relies on the unit-diagonal constraint. It would be valuable to treat rank-deficient affine operators directly, for example through a pseudoinverse or a reduced constraint representation, and to extend the framework to more general affine constraints, and multiple or chordally decomposed PSD blocks.

    \item \emph{Large-scale implementation.}
    For very large problems, explicitly forming or inversing $\cA\cA^*$ and applying the projector may become a bottleneck. Matrix-free projector evaluations, iterative linear solvers with effective preconditioners, parallel treatment of multiple blocks, and warm-start strategies are therefore important implementation directions. More systematic adaptive rules for the penalty parameter, subproblem tolerance, and factorization sizes may further improve robustness.

\end{itemize}

\section*{Declarations}
\subsection*{Funding}
This work was supported by the National Key R\&D Program of China under grant No.~2022YFA1005102 and the Natural Science Foundation of China under grant No.~12571333.

\subsection*{Acknowledgements}
The authors thank the anonymous referees for comments and suggestions that substantially improved the paper.
The authors acknowledge the use of AI for assistance with brainstorming ideas, mathematical development, and drafting the manuscript. The final content, analysis, and conclusions remain the sole responsibility of the
authors.

\subsection*{Competing interests}
The authors have no competing interests to declare that are relevant to the content of this article.

\subsection*{Data availability}
The authors confirm that all data generated or analysed during this study are included in this article.

\bibliographystyle{siamplain}
\bibliography{refer}
\end{document}